\documentclass[preprint,3p,12pt]{elsarticle}

\usepackage{framed,multirow}

\usepackage{amssymb}
\usepackage{latexsym}
\usepackage{graphicx}
\usepackage{caption}

\usepackage{url}

\usepackage{bm}
\usepackage{tabu}
\usepackage{subfigure}
\usepackage{amsmath}
\usepackage{booktabs}
\usepackage[ruled,linesnumbered]{algorithm2e}

\makeatletter
\def\ps@pprintTitle{%
 \let\@oddhead\@empty
 \let\@evenhead\@empty
 \let\@oddfoot\@empty
 \let\@evenfoot\@oddfoot
}
\makeatother

\begin{document}
\captionsetup[figure]{labelfont={bf},labelformat={default},labelsep=period,name={Fig.}}


\begin{frontmatter}

\title{AONN-2: An adjoint-oriented neural network method for PDE-constrained shape optimization}

\author[1]{Xili Wang\fnref{label2}}
\ead{xiliwang@stu.pku.edu.cn}

\author[1]{Pengfei Yin\fnref{label2}}
\ead{pengfeiyin@pku.edu.cn}

\author[1,2]{Bo Zhang}
\ead{zhang_bo@pku.edu.cn}

\author[1,2,3]{Chao Yang\corref{cor1}}
\cortext[cor1]{Corresponding author.}
\ead{chao_yang@pku.edu.cn}

\fntext[label2]{Xili Wang and Pengfei Yin contributed equally to this work.}

\address[1]{School of Mathematical Sciences, Peking University, Beijing 100871, China}
\address[2]{National Engineering Laboratory for Big Data Analysis and Applications, Peking University, Beijing 100871, China}
\address[3]{PKU-Changsha Institute for Computing and Digital Economy, Hunan 410006, China}

\begin{abstract}

Shape optimization has been playing an important role in a large variety of engineering applications. Existing shape optimization methods are generally mesh-dependent and therefore encounter challenges due to mesh deformation. To overcome this limitation, we present a new adjoint-oriented neural network method, AONN-2, for PDE-constrained shape optimization problems. This method extends the capabilities of the original AONN method \cite{yin2023aonn}, which is developed for efficiently solving parametric optimal control problems. AONN-2 inherits the direct-adjoint looping (DAL) framework for computing the extremum of an objective functional and the neural network methods for solving complicated PDEs from AONN. Furthermore, AONN-2 expands the application scope to shape optimization by taking advantage of the shape derivatives to optimize the shape represented by discrete boundary points. AONN-2 is a fully mesh-free shape optimization approach, naturally sidestepping issues related to mesh deformation, with no need for maintaining mesh quality and additional mesh corrections. A series of experimental results are presented, highlighting the flexibility, robustness, and accuracy of AONN-2.

\end{abstract}


\end{frontmatter}

\section{Introduction}


Shape optimization has been playing an important role in a large variety of engineering applications, particularly involving in the design of, e.g., aircraft wings \cite{schmidt2013three}, high-speed train \cite{sun2021aerodynamic}, high-rise building \cite{li2021knowledge}, bridge structure \cite{norato2004geometry}, acoustic devices \cite{Udawalpola2008OptimizationOA}, heat sink \cite{park2004numerical} and biomedical devices \cite{williams2021shape}. Many such problems can be formulated as the minimization of an objective functional defined over a class of admissible domains, which are usually subject to constraints imposed by partial differential equations (PDEs). These PDE-constrained shape optimization problems are notoriously difficult to solve because the dependence of functional on the domain is usually nonconvex. Additionally, even numerically solving the governing PDEs is challenging because the computational domain itself is the unknown variable in the shape optimization problem.

A rather straightforward approach for shape optimization is the heuristic method, which regards a shape optimization problem as a general finite-dimensional optimization problem by parameterizing the shape information into discrete optimization variables. In this way, empirical search algorithms such as the genetic algorithm \cite{akram2021cfd}, the swarm algorithm \cite{ray2004swarm}, and the simulated annealing algorithm \cite{wang2001aerodynamic} can be applied to search the parameterized optimal shape. Heuristic methods are generally easy to implement and good at finding the global optima \cite{liao2021multi}. 
However, they usually require a large number of functional evaluations of different shapes, and are therefore computationally expensive, especially for large scale systems \cite{schmidt2010efficient}. The adjoint methods \cite{kenway2019effective, schulz2016computational, etling2020first, paganini2018higher} take advantage of the gradient information of the objective functional to efficiently update the computational domain for finding the optimal shape. This can be done either on the continuous level or after discretization. The adjoint methods are often combined with the moving mesh approach \cite{paganini2018higher} or the level-set model \cite{allaire2002level}, and have been integrated into several shape optimization toolboxes such as cashocs \cite{blauth2021cashocs} and Fireshape \cite{paganini2021fireshape}. However, the dependency between the mesh and the shape could severely restricts the flexibility of the shape deformation and the update of the level-set function, consequently affecting the accuracy of approximating the optimal shape.
On top of heuristic and adjoint methods, surrogate models can be introduced to alleviate the burden of on-line numerical simulations. A surrogate model is usually constructed off-line so as to directly map design variables to optimization indicators by using, e.g., kriging \cite{sun2021aerodynamic} or polynomial response surface \cite{vidanovic2017aerodynamic}. Still, the construction of a surrogate model typically requires a substantial amount of high-fidelity data obtained through expensive numerical simulations of the governing PDEs.


In recent years, deep learning technologies have gradually been introduced to solve shape optimization problems. For example, fully connected neural networks \cite{zhang2021multi, renganathan2021enhanced, bouhlel2020scalable, xu2021machine}, convolutional neural networks \cite{liao2021multi, xu2021convolutional}, or deep neural operators \cite{shukla2023deep} can be applied to improve the approximation accuracy of surrogate models. And generative adversarial networks \cite{2014gan} can be used to generate more desirable shape samples \cite{zhao2022multi, wang2021airfoil, wu2022missile, li2020efficient}, which could help improve the convergence rate or even the optimization results. However, establishing the generative model needs domain-specific data set, which could severely limit the scope of its application. Moreover, the shape optimization problem can be solved through reinforcement learning such as Q-learning \cite{lampton2010reinforcement}, proximal policy optimization \cite{viquerat2021direct}, or deep deterministic policy gradient \cite{li2021knowledge, yan2019aerodynamic, dai2022aerodynamic} approaches, in which the reward function relies on the shape and guides the policy to provide an action (e.g. design variation) to update the shape. Nevertheless, the reinforcement learning methods often suffer from the meshing failure during the interaction between the agent and the environment \cite{viquerat2021direct}, and they usually can only use several parameters to define the shape, leading to very small shape search space.

In summary, existing methods for shape optimization are usually mesh-dependent, which often suffer from mesh deformation during the optimization process. The mesh deformation could have significant impact on the mesh quality and require additional mesh correction step \cite{biancolini2012mesh}. In this paper, to sidestepping
issues related to mesh deformation, we propose a new adjoint-oriented neural network method, namely AONN-2, for PDE-constrained shape optimization. The AONN method \cite{yin2023aonn} is designed for solving parametric optimal control problems and is primarily applicable to problems featuring explicit control functions within a fixed domain. 
And AONN-2 inherits all the advantages of AONN, including the classical direct-adjoint looping (DAL) framework for computing the extremum of an objective functional \cite{articledaladd} and the emerging neural network methods for solving complicated PDEs \cite{yu2018deep, raissi2019physics, tang2023pinns, sheng2021pfnn, CiCP-32-980}. On top of that, AONN-2  extends the capabilities of AONN by leveraging the shape derivative to optimize the domain's shape, which is represented by a set of discrete boundary points known as shape representing points. In this way, AONN-2 can be applied as a completely mesh-free shape optimization method that is able to naturally avoid issues caused by mesh dependency.

The remaining content of this paper is organized as follows. In Section~\ref{sec:theory}, the basic theory of PDE-constrained shape optimization is introduced. The proposed AONN-2 algorithm is presented in Section~\ref{sec:method}. Then in Section~\ref{sec:result}, a series of numerical experiments are provided to demonstrate the flexibility, robustness, and accuracy of AONN-2. The paper is concluded in Section~\ref{sec:conc}.

\section{Basic theory of PDE-constrained shape optimization}
\label{sec:theory}


PDE-constrained shape optimization can be seen as a special class of optimal control problems, and the corresponding control space is a set of shapes. Let $\Omega\subset \mathbb{R}^{d}$ be a bounded, connected spatial domain with Lipschitz continuous boundary $\partial \Omega$, and $\mathbf{x}\in \Omega$ be the spatial coordinates. A PDE-constrained shape optimization problem can be formulated as:
\begin{equation}
\label{eq:SOP}
\left \{
    \begin{aligned}
      &\min\limits_{(y,\Omega)\in \mathcal{Y}\times \mathcal{U}} J(y,\Omega), \\
      &\text{subject to}\quad \mathbf{F}(y, \Omega)=\mathbf{0},\\
    \end{aligned}
  \right.
\end{equation}
where $J:\mathcal{Y}\times \mathcal{U}\longmapsto \mathbb{R}$ is an objective functional, $\mathcal{Y}$ and $\mathcal{U}$ are the state space and shape space, respectively, and ${y}$ is the solution of the state equation $\mathbf{F}(y,\Omega)=\mathbf{0}$ defined on ${\Omega}$. The state equation is specifically expressed as:




\begin{equation}
\label{eq:state}
\left \{
\begin{aligned}
& \mathbf{F}_I(y, \Omega)(\mathbf{x})=\mathbf{0} \quad \forall \mathbf{x} \in \Omega, \\
& \mathbf{F}_B(y, \Omega)(\mathbf{x})=\mathbf{0} \quad \forall \mathbf{x} \in \partial \Omega.
\end{aligned}
\right.
\end{equation}

\begin{figure}[!htb]
  \centering
  {\includegraphics[width=0.5\textwidth]{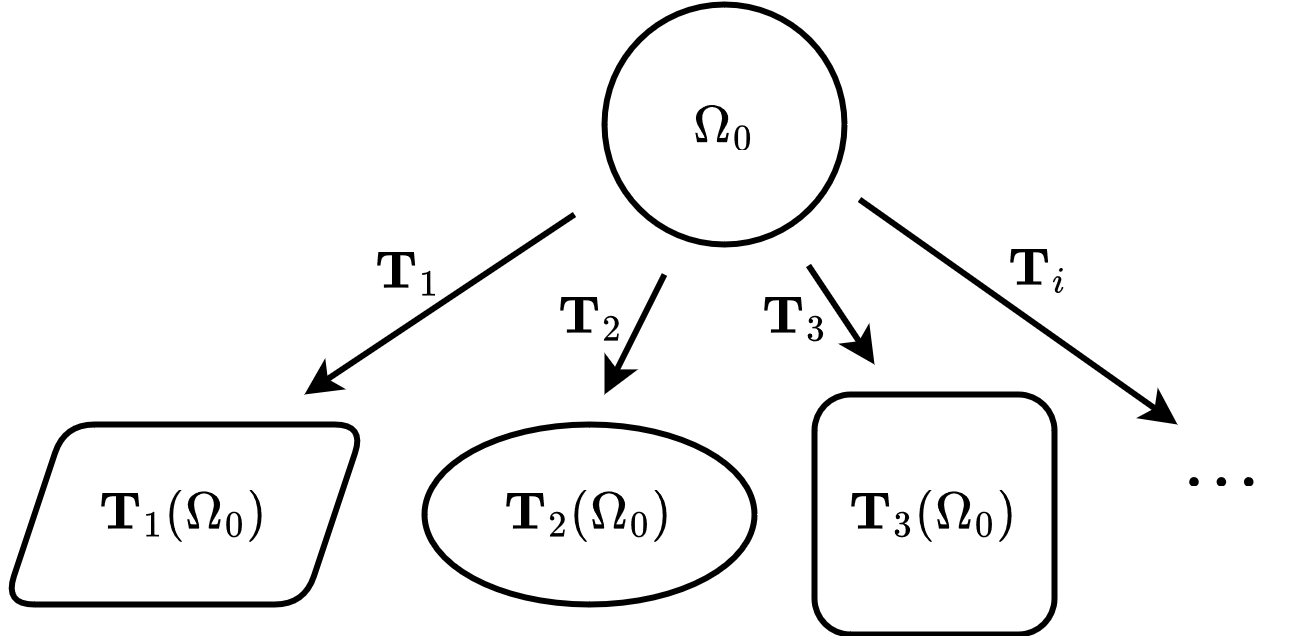}}
  \caption{A set of diffeomorphic geometric transformations $\mathbf{T}_{i}$ (black arrows) applied to an initial domain $\Omega_{0}$.}
  \label{fig:diffeomorphisms}
\end{figure}

In shape optimization, modeling of the shape space is necessary. There are several approaches such as level set method and phase field method. In this work, the shape space $\mathcal{U}$ is modeled as a set of images of diffeomorphic geometric transformations $\mathbf{T}_{i}:\mathbb{R}^{d}\rightarrow\mathbb{R}^{d}$ applied to an initial domain $\Omega_{0}\subset \mathbb{R}^{d}$ \cite{paganini2021fireshape}, that is  
\begin{equation}
\label{eq:shape model}
\mathcal{U}:=\{\Omega=\mathbf{T}_{i}(\Omega_{0}):\mathbf{T}_{i} \in \mathcal{T}_{ad}\}.
\end{equation}
The modeling of the shape space $\mathcal{U}$ is illustrated in Fig.~\ref{fig:diffeomorphisms}. In this way, the domain's boundary can be explicitly described through $\partial (\mathbf{T}_{i}(\Omega_{0}))=\mathbf{T}_{i}(\partial \Omega_{0})$.


For solving problem \eqref{eq:SOP}, the framework of direct-adjoint looping (DAL) is widely adopted \cite{jameson1988aerodynamic, eggl2018gradient, eggl2020shape}. In the process of DAL, an initial domain $\Omega_{0}$ is given firstly, and then the state equation and adjoint equation are solved based on the current domain to get the state variable and adjoint variable. After that, the shape derivative of functional $J$ is computed based on the state and adjoint variables to construct descent direction of the functional. At last, the shape of the domain is updated according to the descent direction. By repeating the above process, the optimized shape can be obtained. In the DAL framework, the computation of the shape derivative and the construction of the descent direction are the two key components.


The shape derivative of functional $J$ in the direction of a vector field $\mathbf{V}$ is defined by:

\begin{equation}
\label{eq:shape derivative}
\mathrm{d}_{\Omega}J(y, \Omega;\mathbf{V}):=\lim_{t\searrow 0}\frac{J(y_{t},\Omega_{t})-J(y,\Omega)}{t},
\end{equation}
where $\Omega_{t}=\{\mathbf{x}+t\mathbf{V}(\mathbf{x}):\mathbf{x}\in\Omega\}$ and $y_{t}$ is the solution to $\mathbf{F}(y_{t},\Omega_{t})=\mathbf{0}$. For efficiently computing the shape derivative, the adjoint equation is usually introduced through constructing the Lagrange functional as: 
\begin{equation}
\label{eq:Lagrange functional}
\mathcal{L}(y,\Omega,p):=J(y,\Omega)+\left(p, \mathbf{F}(y, \Omega)\right)_{L^2(\Omega)},
\end{equation}
and $\mathcal{L}_{y}(y,\Omega,p)=0$ is the adjoint equation, that is 


\begin{equation}
\label{eq:adjoint}
    \begin{aligned}
      \mathbf{F}_{y}(y, \Omega)^*p&=-{J_y}(y,\Omega) \\
    \end{aligned},
\end{equation}
where $p$ is the adjoint variable, which is also known as the Lagrange multiplier. With the state variable $y$ and adjoint variable $p$, which are acquired by solving the state equation and adjoint equation, respectively, the shape derivative of functional $J$ in direction $\mathbf{V}$ can be computed by: 
\begin{equation}
\label{eq:comput shape derivative}
\mathrm{d}_{\Omega}J(y,\Omega;\mathbf{V})=\mathcal{L}_{\Omega}(y,\Omega, p;\mathbf{V}).
\end{equation}

After obtaining the shape derivative, the shape gradient $\nabla\hat{J}(\Omega): \partial\Omega \rightarrow \mathbb{R}$, which satisfies $\mathrm{d}_{\Omega}J(y,\Omega;\mathbf{V}) = (\nabla\hat{J}(\Omega),\mathbf{V}\cdot\mathbf{n})_{L^2(\partial\Omega)}$, can be naturally employed to construct the descent direction as $-\nabla\hat{J}(\Omega)\mathbf{n}$ on $\partial\Omega$, where $\mathbf{n}$ is the unit outward normal vector to $\Omega$, and $\hat{J}(\Omega) = J(y,\Omega)$ is the reduced form of the objective functional. Unfortunately, as demonstrated in Section~\ref{sec:method}, this descent direction often exhibits low regularity sometimes. To enhance the regularity of the descent direction, the classical Hilbertian method \cite{allaire2004structural, allaire2021shape} is frequently adopted. In this method, the regularization equation \eqref{eq:regularization} needs to be solved to get the regularized descent direction $\mathbf{\Phi}$.
\begin{equation}
\label{eq:regularization}
\left(\mathbf{\Phi},\mathbf{V}\right)_{[H^{1}(\Omega)]^d}=-\mathrm{d}_{\Omega}J(y,\Omega;\mathbf{V}) \quad \forall \; \mathbf{V}\in W^{1,\infty}(\mathbb{R}^{d},\mathbb{R}^{d}).
\end{equation}
Different forms of inner products can be employed in equation \eqref{eq:regularization}, and the $H^1$ inner product is used in this work. The strong form of equation \eqref{eq:regularization} is denoted by $\mathbf{G}(y,p,\mathbf{\Phi},\Omega) = \mathbf{0}$ with Neumann boundary condition, which is specifically expressed as: 

\begin{equation}
\label{reg_eq}
\left \{
\begin{aligned}
\mathbf{G}_I(\mathbf{\Phi},\Omega) & := - \Delta \mathbf{\Phi}+\mathbf{\Phi} &=\mathbf{0} & & \text { in } \Omega,\\
\mathbf{G}_B(y,p,\mathbf{\Phi},\Omega) & := \partial_{\mathbf{n}} \mathbf{\Phi} + \nabla\hat{J}(\Omega) \mathbf{n}&=\mathbf{0}& & \text { on } \partial \Omega. 
\end{aligned} 
\right.
\end{equation}

Finally, the domain $\Omega$ can be updated by using $\mathbf{\Phi}$ as the descent direction of objective functional. In the following section, AONN-2 will be introduced, which is developed based on the aforementioned theory.

\section{Methodology}
\label{sec:method}

The adjoint-oriented neural network method (AONN) \cite{yin2023aonn} was proposed to solve parametric optimal control problems without directly solving the complex Karush–Kuhn–Tucker system with various penalty terms. It combines the DAL framework with neural network approaches for solving PDEs. In AONN, the control variable to be optimized is restricted to the form of a function, and without considering parameters, the control variable takes the form of a function $u(\mathbf{x})$ that is defined on a fixed computational domain $\Omega$ throughout the optimization process. Specifically, in AONN the solution of the optimal system is approximated by iteratively updating the neural networks related to the state variable $y$, adjoint variable $p$ and control variable $u$. In the context of shape optimization, it becomes necessary to consider how to parameterize and update the shape of the domain, as well as how to solve the PDEs in the changing computational domain. Addressing these challenges is the crucial target of this work.


\begin{figure}[!htb]
  \centering
  \subfigure[]
  {\includegraphics[width=0.28\textwidth]{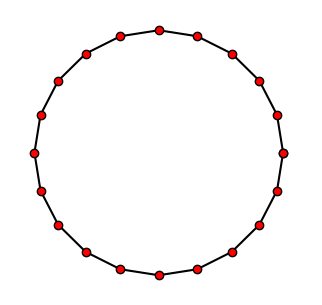}}
  \subfigure[]
  {\includegraphics[width=0.28\textwidth]{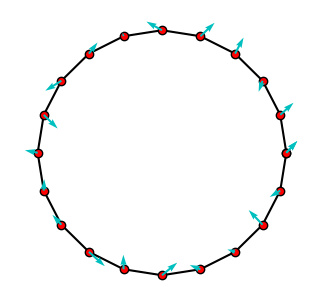}}
  \subfigure[]
  {\includegraphics[width=0.28\textwidth]{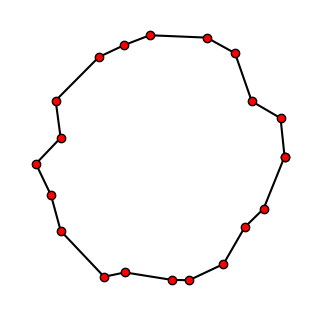}}
  \caption{Representing and updating shapes. (a) An initial shape is represented by discrete boundary points (i.e. shape representing points). (b) The descent direction $\mathbf{\Phi}$ (blue arrows) is used to guide the moving of shape representing points. (c) The updated shape is attained with retaining the connected relation about the shape representing points.}
  \label{fig:bpm}
\end{figure}


We extend AONN to AONN-2 that can solve PDE-constrained shape optimization problems. Firstly, a shape is represented by discrete points $\mathbf{x}$ on the shape's boundary, which are named as shape representing points in this paper. Based on this representation, the shape can be updated according to the descent direction $\mathbf{\Phi}$ as $\mathbf{x}^{\ast} = \mathbf{x}+\alpha\mathbf{\Phi}(\mathbf{x})$, where $\alpha$ is the step length. The schematic diagram about representing and updating the shape is shown in Fig.~\ref{fig:bpm}. To obtain the descent direction $\mathbf{\Phi}$, one needs to solve three PDEs: the state equation \eqref{eq:state}, the adjoint equation \eqref{eq:adjoint} and the regularization equation \eqref{reg_eq} as exhibited in Section~\ref{sec:theory}. In AONN-2, these PDEs are solved based on physics-informed neural networks (PINNs) \cite{raissi2019physics}, which are intrinsically mesh-free. After updating the shape, it only needs to resample the collocation points in the updated domain for solving the three PDEs in the next step. However, in traditional mesh-dependent shape optimization methods, it is necessary to generate internal meshes and deform them, and it usually needs to correct the meshes to guarantee their quality \cite{etling2020first, iglesias2018two}. The whole process of mesh processing is complex and time-consuming. The proposed AONN-2 avoids the mesh processing, and focuses on the movement of the shape representing points. The cost of resampling collocation points is much lower than the mesh processing, and the movement of shape representing points is also more flexible than mesh deformation. 

For solving the state equation, adjoint equation and regularization equation based on PINNs, the corresponding neural networks are established as $\tilde{y}(\mathbf{x};\boldsymbol{\theta}_{y})$, $\tilde{p}(\mathbf{x};\boldsymbol{\theta}_{p})$ and $\tilde{\mathbf{\Phi}}(\mathbf{x};\boldsymbol{\theta}_{\mathbf{\Phi}})$ with network parameters $\boldsymbol{\theta}_{y}$, $\boldsymbol{\theta}_{p}$ and $\boldsymbol{\theta}_{\mathbf{\Phi}}$ respectively. The related three loss functions are defined as:

\begin{equation}
\label{eq:loss state}
L_{s}(\boldsymbol{\theta}_{y},\Omega) = \left(\frac{1}{N}\sum_{i=1}^{N}|r_{s_I}(\tilde{y}(\mathbf{x}_{I}^{i};\boldsymbol{\theta}_{y}),\Omega)|^{2} + \frac{\lambda_s}{M}\sum_{i=1}^{M}|r_{s_B}(\tilde{y}(\mathbf{x}_{B}^{i};\boldsymbol{\theta}_{y}),\Omega)|^{2}\right)^{\frac{1}{2}},
\end{equation}

\begin{equation}
\label{eq:loss adjoint}
\begin{aligned}
L_{a}(\boldsymbol{\theta}_{y},\boldsymbol{\theta}_{p},\Omega) = & \left(\frac{1}{N}\sum_{i=1}^{N}|r_{a_I}(\tilde{y}(\mathbf{x}_{I}^{i};\boldsymbol{\theta}_{y}),\tilde{p}(\mathbf{x}_{I}^{i};\boldsymbol{\theta}_{p}),\Omega)|^{2} + \right.\\
& \left. \frac{\lambda_a}{M}\sum_{i=1}^{M}|r_{a_B}(\tilde{y}(\mathbf{x}_{B}^{i};\boldsymbol{\theta}_{y}),\tilde{p}(\mathbf{x}_{B}^{i};\boldsymbol{\theta}_{p}),\Omega)|^{2}\right)^{\frac{1}{2}},
\end{aligned}
\end{equation}

\begin{equation}
\label{eq:loss regularization}
\begin{aligned}
L_{r}(\boldsymbol{\theta}_{y},\boldsymbol{\theta}_{p},\boldsymbol{\theta}_{\mathbf{\Phi}},\Omega) = &\left(\frac{1}{N}\sum_{i=1}^{N}|r_{r_I}(\tilde{y}(\mathbf{x}^{i}_{I};\boldsymbol{\theta}_{y}),\tilde{p}(\mathbf{x}^{i}_{I};\boldsymbol{\theta}_{p}), \tilde{\mathbf{\Phi}}(\mathbf{x}^{i}_{I};\boldsymbol{\theta}_{\mathbf{\Phi}}),\Omega)|^{2} +\right.\\
&\left. \frac{\lambda_r}{M}\sum_{i=1}^{M}|r_{r_B}(\tilde{y}(\mathbf{x}^{i}_{B};\boldsymbol{\theta}_{y}),\tilde{p}(\mathbf{x}^{i}_{B};\boldsymbol{\theta}_{p}), \tilde{\mathbf{\Phi}}(\mathbf{x}^{i}_{B};\boldsymbol{\theta}_{\mathbf{\Phi}}),\Omega)|^{2} \right)^{\frac{1}{2}},
\end{aligned}
\end{equation}
where $\{\mathbf{x}^{i}_{I}\}_{i=1}^{N}$ denote $N$ collocation points in the domain $\Omega$, $\{\mathbf{x}^{i}_{B}\}_{i=1}^{M}$ denote $M$ collocation points on the domain's boundary $\partial\Omega$, which can employ the shape representing points, and $r_{s}$, $r_{a}$ and $r_{r}$ represent the residuals for the state equation, adjoint equation and regularization equation as:
\begin{subequations}
\begin{align}
r_{s}(y,\Omega)&=\mathbf{F}(y,\Omega), \label{eq:residul state}\\
r_{a}(y,p,\Omega)&=\mathbf{F}_{y}(y,\Omega)^{\ast}p + J_{y}(y,\Omega), \label{eq:residul adjoint}\\
r_{r}(y,p,\mathbf{\Phi},\Omega)&=\mathbf{G}(y,p,\mathbf{\Phi},\Omega).\label{eq:residul regularization}
\end{align}
\end{subequations}
According to equation \eqref{eq:state}, the state residual $r_{s}$ \eqref{eq:residul state} is divided to the interior residual $r_{s_I}$ and boundary residual $r_{s_B}$. Similarly, the adjoint residual $r_{a}$ \eqref{eq:residul adjoint} is divided to $r_{a_I}$ and $r_{a_B}$, and the regularization residual $r_{r}$ \eqref{eq:residul regularization} is divided to $r_{r_I}$ and $r_{r_B}$. The corresponding weights for the boundary residuals are denoted as $\lambda_s$, $\lambda_a$, and $\lambda_r$, respectively. The updating of the neural network parameters $\boldsymbol{\theta}_{y},\boldsymbol{\theta}_{p},\boldsymbol{\theta}_{\mathbf{\Phi}}$ can be realized by automatic differentiation and back propagation with deep learning libraries, such as PyTorch \cite{paszke2017automatic} and TensorFlow \cite{abadi2016tensorflow}. 

\begin{figure}[!htb]
  \centering
  \subfigure[]
  {\includegraphics[width=0.275\textwidth]{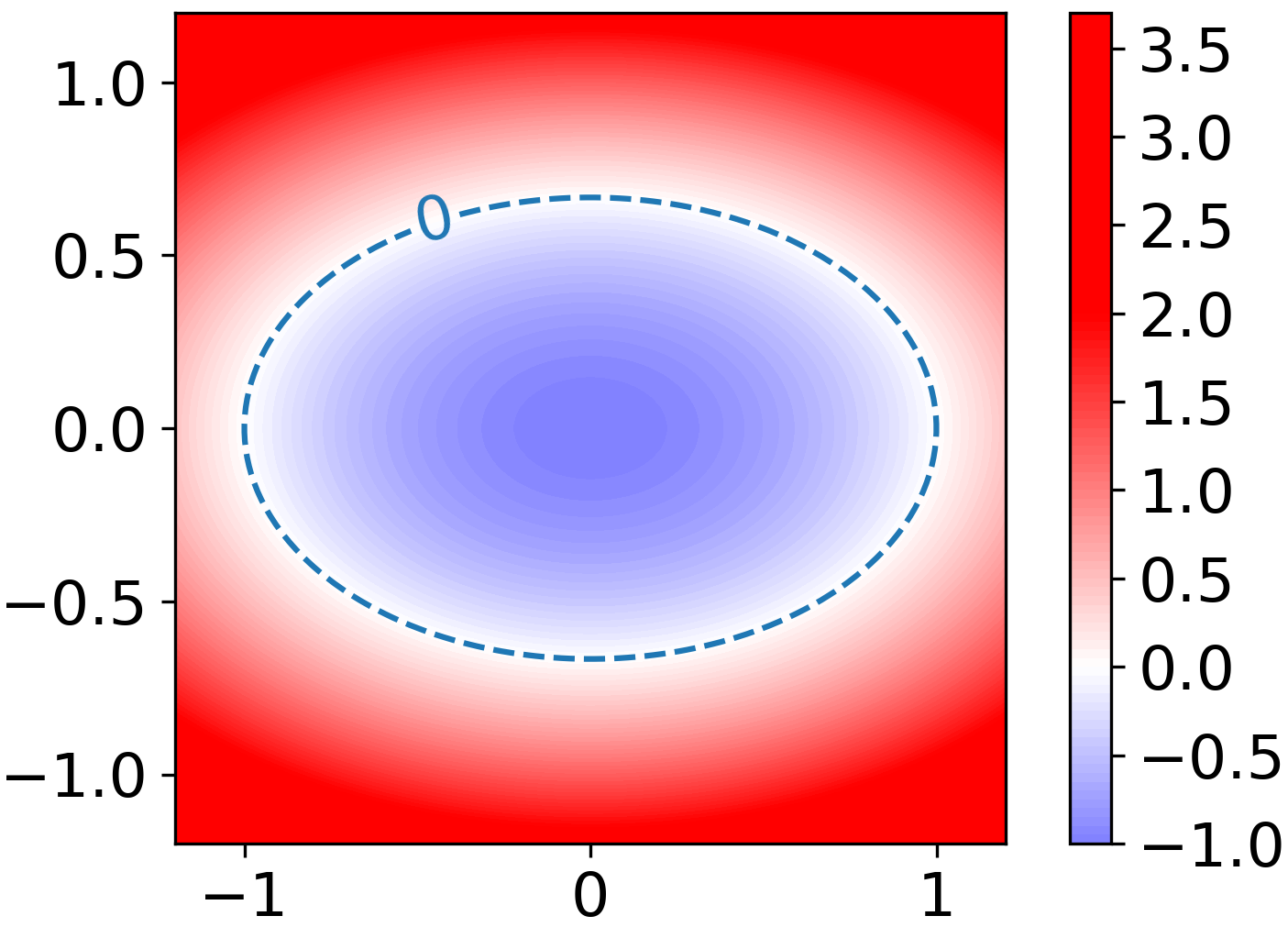}\label{fig:example-a}}
  \subfigure[]
  {\includegraphics[width=0.22\textwidth]{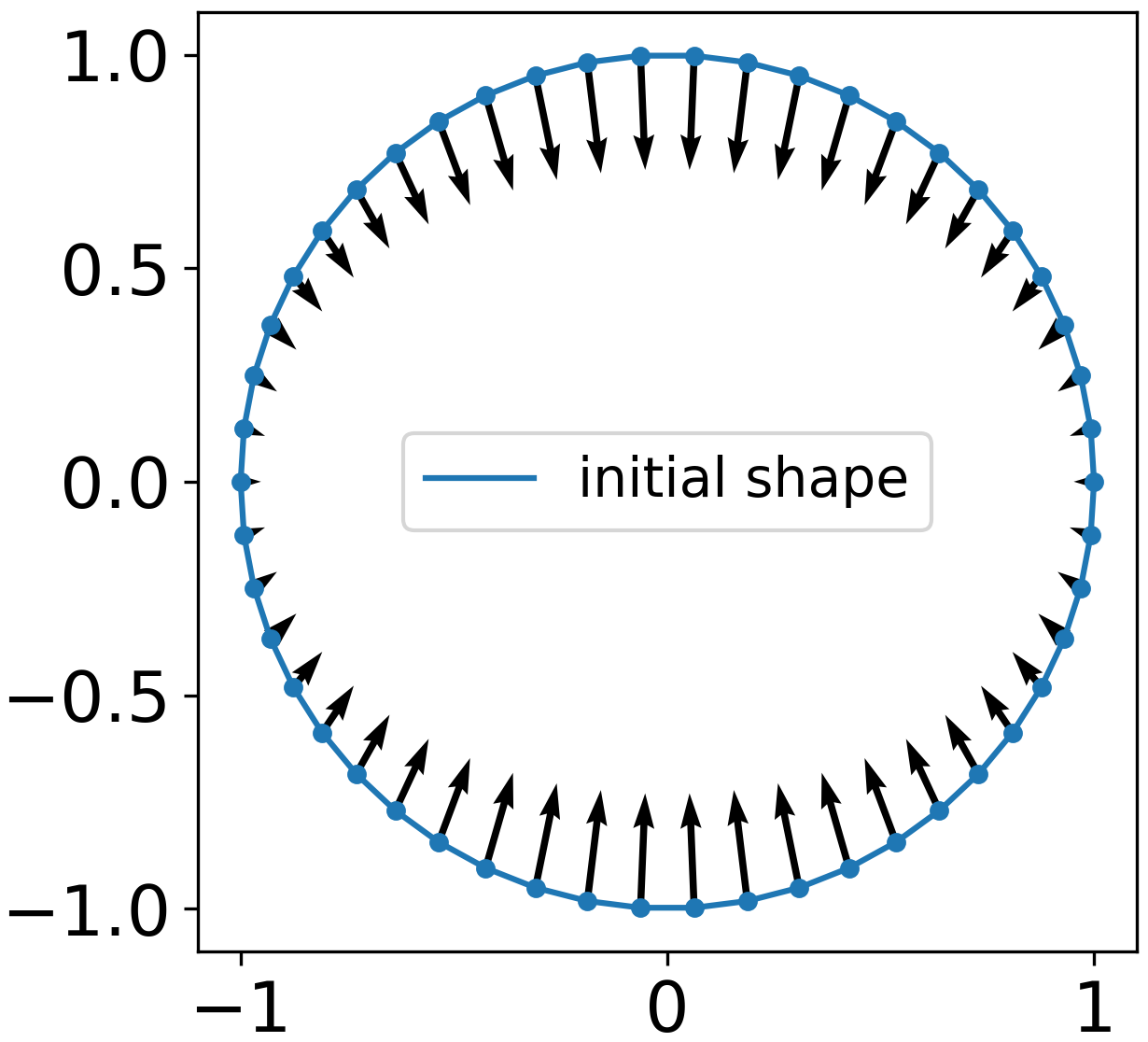}\label{fig:example-b}}
  \subfigure[]
  {\includegraphics[width=0.22\textwidth]{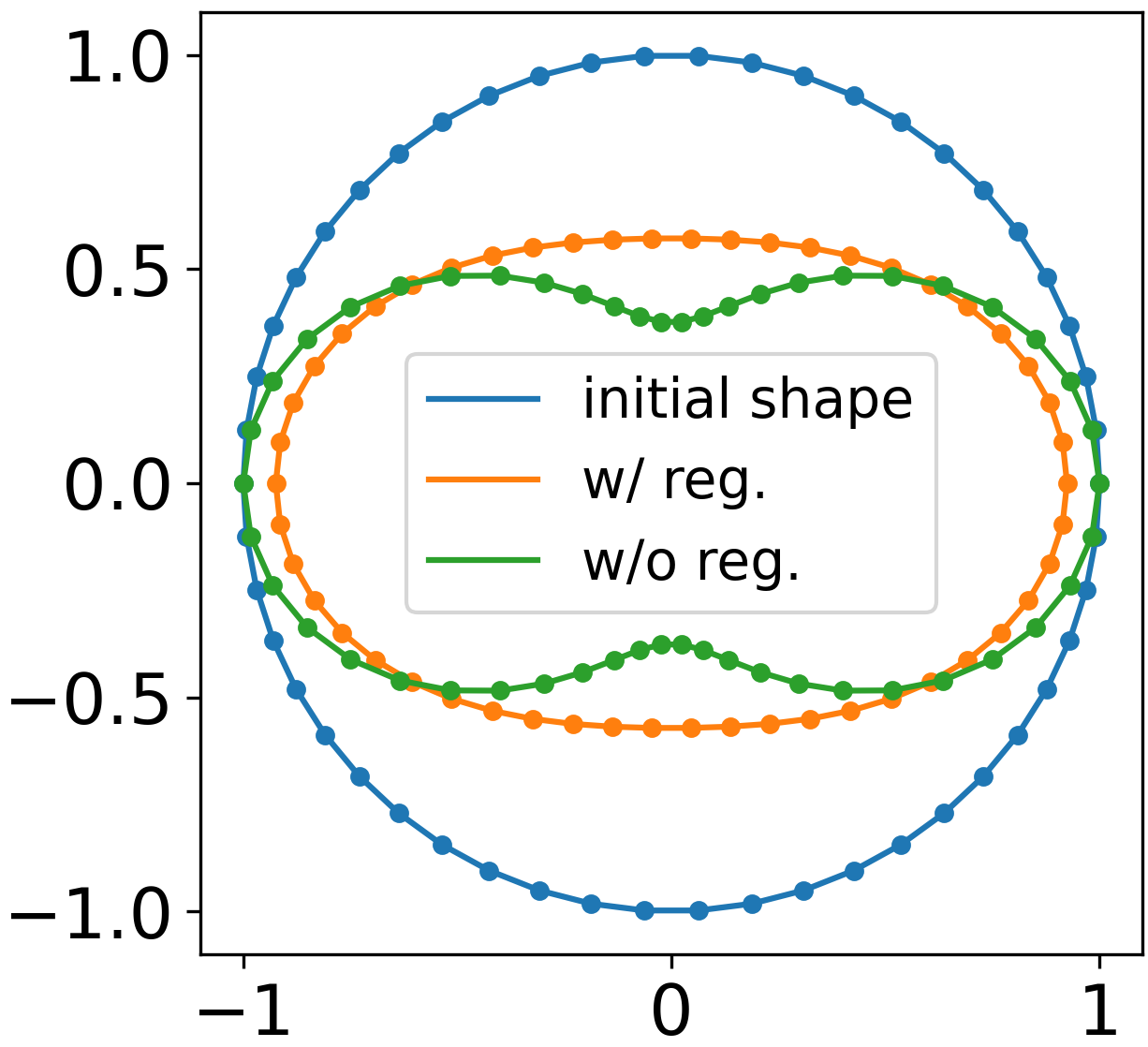}\label{fig:example-c}}
  \subfigure[]
  {\includegraphics[width=0.22\textwidth]{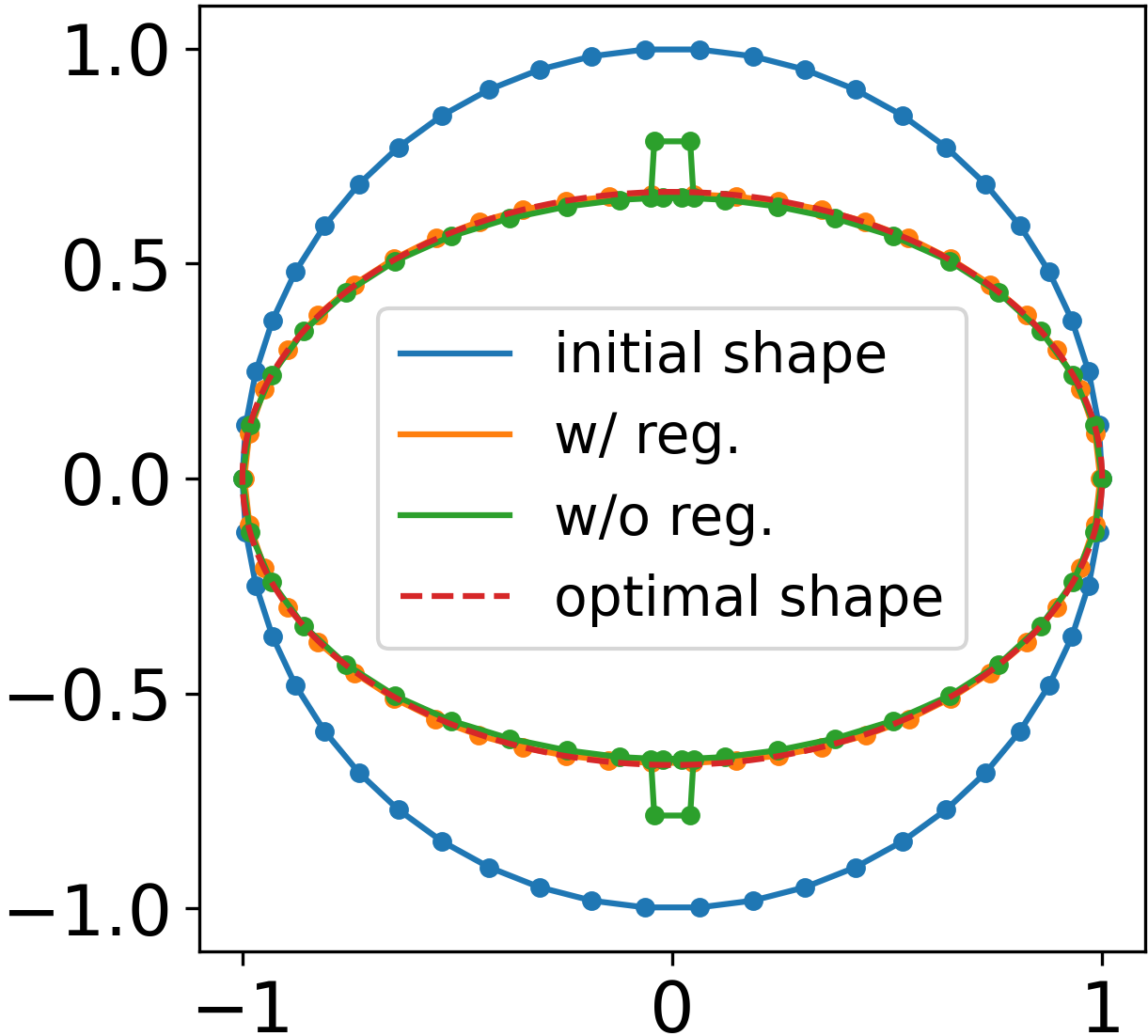}\label{fig:example-d}}
  \caption{Comparison of shape updating with regularization (w/ reg.) and without regularization (w/o reg.). (a) The optimal shape (i.e. the zero level set of $f$), which is shown as the blue dotted line. (b) Deformation of each shape representing point of initial shape (w/ reg.). (c) Shapes after once update. (d) Shapes after five updates.}
  \label{fig:demo}
\end{figure}

In order to verify the necessity of the regularization of the descent direction in AONN-2, a simple two-dimensional unconstrained shape optimization problem is considered:

\begin{equation}
\begin{aligned}
  &\min\limits_{\Omega}J(f,\Omega):=\int_{\Omega}f \,\mathrm{d}\mathbf{x}, \\
\end{aligned}
\end{equation}
where $f\left(x_{1}, x_{2}\right)=x_{1}^{2}+\left(3 x_{2} / 2\right)^{2}-1$ and the optimal $\Omega$ is searched for minimizing the integral value. Obviously, the integral value is minimized when the shape's boundary is set to the zero level set of $f(x_{1},x_{2})$. For regularizing the descent direction, the corresponding regularization equation \eqref{demo:reg_eq} is solved:
\begin{equation}
\label{demo:reg_eq}
\left\{\begin{aligned}
- \Delta \mathbf{\Phi}+\mathbf{\Phi} &=\mathbf{0} & \text { in } \Omega,\\
\partial_{\mathbf{n}} \mathbf{\Phi} + f \mathbf{n}&=\mathbf{0} & \text { on } \partial \Omega.
\end{aligned}\right.
\end{equation}
In this example, the initial shape is set to an unit circle. The comparison of shape updating with and without regularization is shown in Fig.~\ref{fig:demo}. From this figure, we can see that the shape updating without regularization has poor smoothness and may lead to divergence during the iteration process. Thus it is necessary to improve the regularity of descent direction.

\begin{figure}[!htb]
  \centering
  {\includegraphics[height=8.0cm]{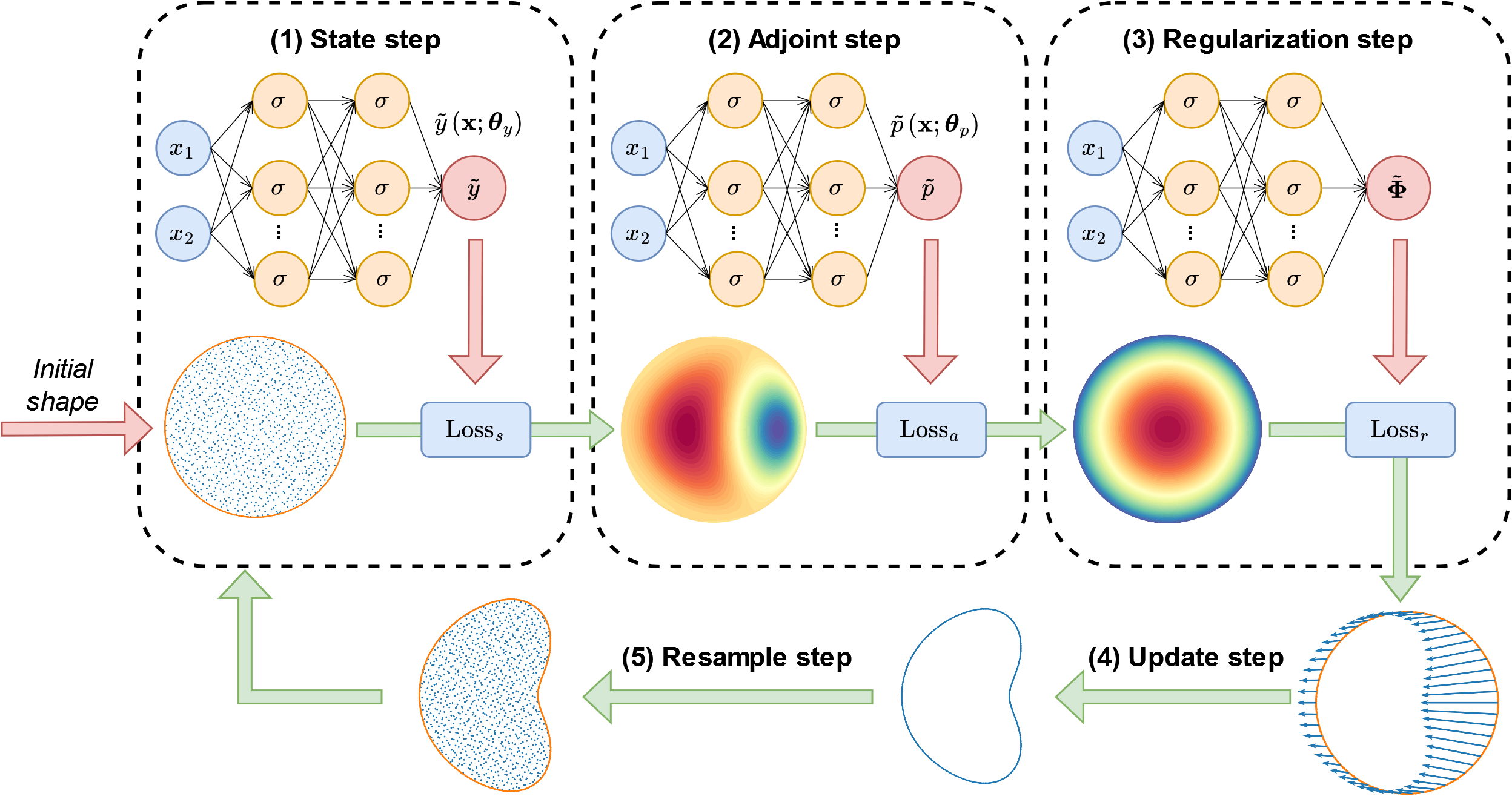}}
  \caption{The schematic diagram about the computation process of AONN-2. It mainly contains state step, adjoint step, regularization step, update step and resample step.}   
  \label{fig:working flow}
\end{figure}

In summary, the computation process of AONN-2 consists of a looping of five steps: training $\tilde{y}(\mathbf{x};\boldsymbol{\theta}_{y})$ (state step), training $\tilde{p}(\mathbf{x};\boldsymbol{\theta}_{p})$ (adjoint step), training $\tilde{\mathbf{\Phi}}(\mathbf{x};\boldsymbol{\theta}_{\mathbf{\Phi}})$ (regularization step), updating the shape representing points (update step), and resampling the collocation points (resample step), which is illustrated in Fig.~\ref{fig:working flow}. More specifically, starting with initial neural networks $\tilde{y}(\mathbf{x};\boldsymbol{\theta}_{y}^{0})$,  $\tilde{p}(\mathbf{x};\boldsymbol{\theta}_{p}^{0})$, $\tilde{\mathbf{\Phi}}(\mathbf{x};\boldsymbol{\theta}_{\mathbf{\Phi}}^{0})$ and an initial domain $\Omega_{0}$, the state equation is solved through minimizing loss function $L_{s}(\boldsymbol{\theta}_{y}, \Omega_{0})$ to get the solution $\tilde{y}(\mathbf{x};\boldsymbol{\theta}_{y}^{1})$. Following this, with $\tilde{y}(\mathbf{x};\boldsymbol{\theta}_{y}^{1})$, the loss function $L_{a}(\boldsymbol{\theta}_{y}^{1},\boldsymbol{\theta}_{p}, \Omega_{0})$ is minimized by solving the adjoint equation to get the solution $\tilde{p}(\mathbf{x};\boldsymbol{\theta}_{p}^{1})$. After that, with the solutions of state and adjoint equations, the descent direction $\tilde{\mathbf{\Phi}}(\mathbf{x};\boldsymbol{\theta}_{\mathbf{\Phi}}^{1})$ of the objective functional can be obtained by solving the regulation equation through minimizing the loss function $L_{r}(\boldsymbol{\theta}_{y}^{1},\boldsymbol{\theta}_{p}^{1},\boldsymbol{\theta}_{\mathbf{\Phi}}, \Omega_{0})$. Then, $\tilde{\mathbf{\Phi}}(\mathbf{x};\boldsymbol{\theta}_{\mathbf{\Phi}}^{1})$ is employed to guide the updating of the shape representing points of $\Omega_{0}$ to form a new domain $\Omega_{1}$. Lastly, the collocation points in $\Omega_{1}$ are resampled. The next iteration begins with the network parameters $\boldsymbol{\theta}_{y}^{1},\boldsymbol{\theta}_{p}^{1},\boldsymbol{\theta}_{\mathbf{\Phi}}^{1}$ and $\Omega_{1}$. The updates of the network parameters and the shape are expressed as follows:

\begin{equation}
\nonumber
\begin{aligned}
\boldsymbol{\theta}_{y}^{k+1}&=\mathop{\arg\min}\limits_{\boldsymbol{\theta}_{y}}L_{s}(\boldsymbol{\theta}_{y},\Omega_{k}),\\
\boldsymbol{\theta}_{p}^{k+1}&=\mathop{\arg\min}\limits_{\boldsymbol{\theta}_{p}}L_{a}(\boldsymbol{\theta}_{y}^{k+1},\boldsymbol{\theta}_{p},\Omega_{k}),\\
\boldsymbol{\theta}_{\mathbf{\Phi}}^{k+1}&=\mathop{\arg\min}\limits_{\boldsymbol{\theta}_{\mathbf{\Phi}}}L_{r}(\boldsymbol{\theta}_{y}^{k+1},\boldsymbol{\theta}_{p}^{k+1},\boldsymbol{\theta}_{\mathbf{\Phi}},\Omega_{k}),\\
\end{aligned}
\end{equation}
and
\begin{equation}
\nonumber
\partial\Omega_{k+1}=\partial\Omega_{k}+\alpha \tilde{\mathbf{\Phi}}(\partial\Omega_{k};\boldsymbol{\theta}_{\mathbf{\Phi}}^{k+1}).
\end{equation}

The whole computation process is presented in Algorithm~\ref{alg3}. It is worth mentioning that the step size $\alpha$ is crucial for the convergence of AONN-2. A large step size may lead to divergence, while a small step size may cause slow convergence. As a compromise, a moderate starting step size $\alpha_{0}$ can be chosen and a decay strategy with decay factor $\gamma$ can be adopted.

\begin{algorithm}[!htb]
	\caption{Computation process of the AONN-2 algorithm.}
	\label{alg3}
	\KwIn{Shape representing points of initial domain $\Omega_{0}$, initial network parameters $\boldsymbol{\theta}_{y}^{0},\boldsymbol{\theta}_{p}^{0},\boldsymbol{\theta}_{\mathbf{\Phi}}^{0}$, initial step size $\alpha_{0}$, decay factor $\gamma\in(0,1)$, number of iterations $K$, $N$ collocation points $\{\mathbf{x}^{i}_{I}\}_{i=1}^{N}$ in $\Omega_{0}$ and $M$ collocation points $\{\mathbf{x}^{i}_{B}\}_{i=1}^{M}$ on $\partial\Omega_0$ which can use the shape representing points.}
	\KwOut{Optimized domain $\Omega_{K}$, optimized state $\tilde{y}(\mathbf{x};\boldsymbol{\theta}_{y}^{K})$.}  
	\BlankLine
	
	\While{$k < K$}{
		$\boldsymbol{\theta}_{y}^{k+1}\leftarrow \mathop{\arg\min}\limits_{\boldsymbol{\theta}_{y}}L_{s}(\boldsymbol{\theta}_{y},\Omega_{k})$: Train the neural network $\tilde{y}(\mathbf{x};\boldsymbol{\theta}_{y})$ from the initialization $\boldsymbol{\theta}_{y}^{k}$.

        $\boldsymbol{\theta}_{p}^{k+1}\leftarrow \mathop{\arg\min}\limits_{\boldsymbol{\theta}_{p}}L_{a}(\boldsymbol{\theta}_{y}^{k+1},\boldsymbol{\theta}_{p},\Omega_{k})$: Train the neural network $\tilde{p}(\mathbf{x};\boldsymbol{\theta}_{p})$ from the initialization $\boldsymbol{\theta}_{p}^{k}$.

        $\boldsymbol{\theta}_{\mathbf{\Phi}}^{k+1}\leftarrow \mathop{\arg\min}\limits_{\boldsymbol{\theta}_{\mathbf{\Phi}}}L_{r}(\boldsymbol{\theta}_{y}^{k+1},\boldsymbol{\theta}_{p}^{k+1},\boldsymbol{\theta}_{\mathbf{\Phi}},\Omega_{k})$: Train the neural network $\tilde{\mathbf{\Phi}}(\mathbf{x};\boldsymbol{\theta}_{\mathbf{\Phi}})$ from the initialization $\boldsymbol{\theta}_{\mathbf{\Phi}}^{k}$.

        $\partial\Omega_{k+1}\leftarrow\partial\Omega_{k}+\alpha_{k} \tilde{\mathbf{\Phi}}(\partial\Omega_{k};\boldsymbol{\theta}_{\mathbf{\Phi}}^{k+1})$: Update the shape representing points 
        
        Resample $N$ collocation points $\{\mathbf{x}^{i}_{I}\}_{i=1}^{N}$ in $\Omega_{k+1}$.
        
        $\alpha_{k+1}\leftarrow \gamma \alpha_{k}$.

        $k\leftarrow k+1$.
	}
	
\end{algorithm}

\section{Numerical experiments}
\label{sec:result}
In this section, a series of numerical experiments are carried out to demonstrate the effectiveness of AONN-2. In these experiments, different types of PDE constraints and objective functionals are considered, and AONN-2 are compared with the shape optimization toolbox Fireshape \cite{paganini2021fireshape}, which is based on the adjoint method with finite element discretization. In AONN-2, ResNet \cite{he2016deep} with sinusoid activation function is adopted as the neural network model, and each residual block is comprised of two fully connected layers and a residual connection. Unless otherwise specified, the quasi Monte-Carlo method and uniform sampling method are respectively used to generate the collocation points in the domain and on the boundary by calling the SciPy module \cite{virtanen2020scipy}. For training the neural networks, the Broyden-Fletcher-Goldfarb-Shanno (BFGS) algorithm with a strong Wolfe line search strategy is employed based on PyTorch \cite{paszke2017automatic}. All the training is performed on a server equipped with a Geforce RTX 2080 GPU and based on 64-bit floating point precision. In Fireshape, Limited-Memory BFGS (L-BFGS) method within Rapid Optimization Library (ROL) \cite{ridzal2017rapid} is used to perform the updating of shape. The codes accompanying this manuscript will be published in GitHub (\url{https://github.com/SillyWWang/nn4shape}).

\subsection{Model problem constrained by Poisson equation}
\begin{figure}[!htb]
  \centering
   \subfigure[Initialization with circle]
  {\includegraphics[width=0.96\textwidth]{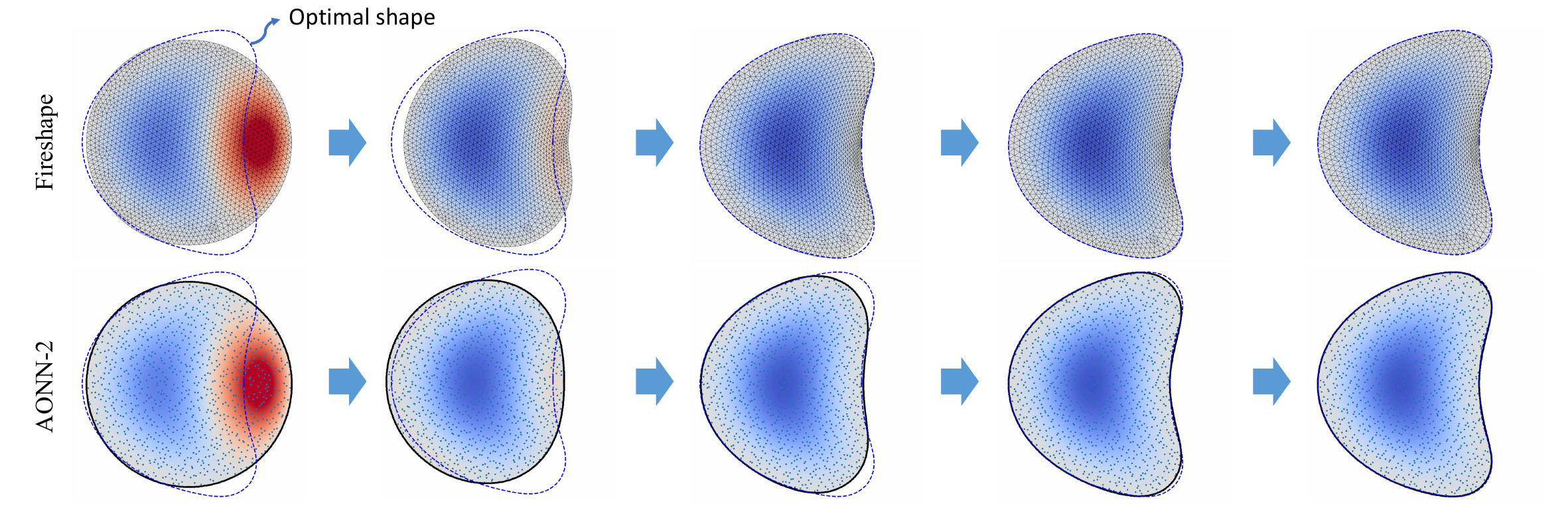}}
   \subfigure[Initialization with ellipse]
  {\includegraphics[width=0.96\textwidth]{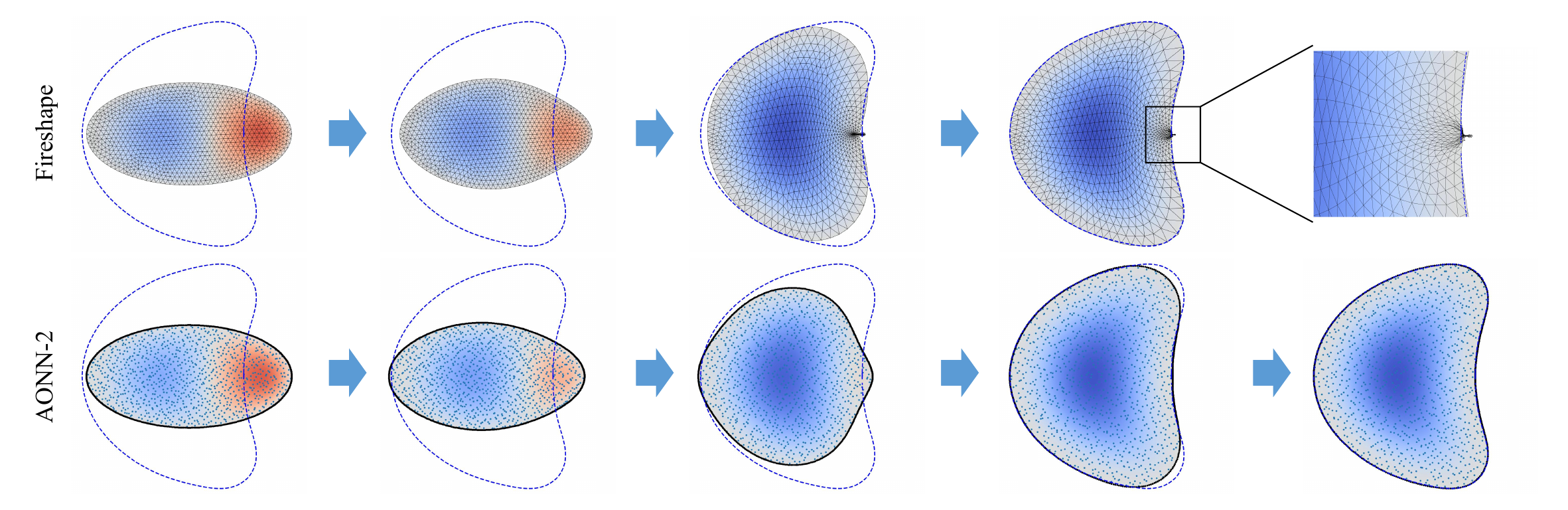}}
  \subfigure[Initialization with rectangle]
  {\includegraphics[width=0.96\textwidth]{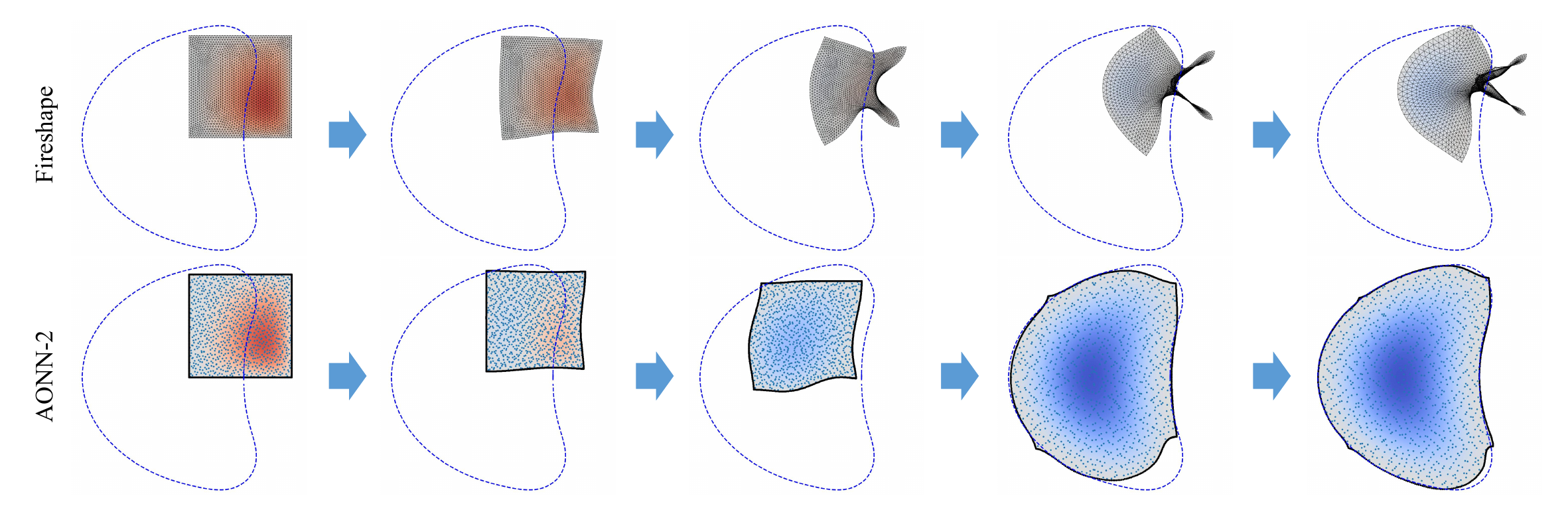}}
  \caption{Shape optimization results by Fireshape and AONN-2, with three different initial shapes: (a) circle (b) ellipse (c) rectangle. The color represents the magnitude of $y$.}
  \label{fig:numerical_example_1}
\end{figure}
We start with a model problem constrained by the Poisson equation with homogeneous Dirichlet boundary condition \cite{etling2020first, blauth2021cashocs}. This problem can be expressed by:
\begin{equation}
\label{eq:Poisson optimization}
\left \{
    \begin{aligned}
      &\min\limits_{y, \Omega}J(y,\Omega):=\int_{\Omega}y \,\mathrm{d}\mathbf{x}, \\
      &\text{subject to} \left \{ 
                   \begin{aligned}
                   -\Delta y &= f && \text{ in } \Omega,\\
                    y &= 0 && \text{ on } \partial\Omega,
                   \end{aligned}
                   \right.\\
    \end{aligned}
  \right.
\end{equation}
where $f(x_{1},x_{2})=2.5(x_{1}+0.4-x_{2}^{2})^{2} + x_{1}^{2}+x_{2}^{2}-1$ according to \cite{blauth2021cashocs}. The shape derivative of functional $J$ is:
\begin{equation}
\label{eq:Poisson_shape_derivate}
\mathrm{d}_{\Omega}J(y,\Omega;\mathbf{V}) = \int_{\partial\Omega}(\partial_\mathbf{n}y\partial_\mathbf{n}p)\mathbf{n}\cdot \mathbf{V} \mathrm{d}s,
\end{equation}
where $\mathbf{n}$ is the unit outward normal vector to $\Omega$ and $p$ is the solution of the corresponding adjoint equation:
\begin{equation}
\left \{
  \begin{aligned}
  -\Delta p &=1 && \text{ in }\Omega,\\
  p&=0 && \text{ on } \partial \Omega.
  \end{aligned}
\right.
\end{equation}
And to get the regularized descent direction $\mathbf{\Phi}$, the following regularization equation needs to be solved:
\begin{equation}
  \left\{
    \begin{aligned}
      -\Delta \mathbf{\Phi}+\mathbf{\Phi}&=\mathbf{0} && \text{ in }\Omega,\\
      \partial_\mathbf{n}\mathbf{\Phi}+(\partial_\mathbf{n}y\partial_\mathbf{n}p)\mathbf{n}&=\mathbf{0} &&\text{ on } \partial \Omega.
    \end{aligned}
  \right.
\end{equation}

In AONN-2, $y, p$ and $\mathbf{\Phi}$ are respectively expressed by three neural networks $\tilde{y}(\mathbf{x};\theta_{y})$, $\tilde{p}(\mathbf{x};\theta_{p})$ and $\tilde{\mathbf{\Phi}}(\mathbf{x};\theta_{\mathbf{\Phi}})$, and each of them have 2 residual blocks with 10 neurons per hidden layer. These networks take spatial coordinates $(x_{1},x_{2})$ as input and output $y,p,\mathbf{\Phi}$ respectively. For solving the three PDEs based on PINNs, 1000 collocation points inside the domain and 500 collocation points on the boundary are sampled, and the boundary collocation points are also used as the shape representing points. In this problem, the initial domain $\Omega_{0}$ is respectively set to circle, ellipse and rectangle as shown in Fig.~\ref{fig:numerical_example_1}. As a comparison, Fireshape is employed for solving this problem. In Fireshape, the P1 Lagrange finite element method is used for spatial discretization and the Krylov subspace method in PETSc is adopted to solve the discretized system. In particular, a total of 3120, 1776, and 3872 triangles are respectively used to discrete the computational domain for the cases with initial shapes of circle, ellipse, and rectangle. The maximal iteration number $K$ in Fireshape and AONN-2 is respectively set to 20 and 50. In order to evaluate the performance of both methods, the optimization result from \cite{blauth2021cashocs} is adopted as a reference.


Starting from different initializations, the results during the whole optimization processes of Fireshape and AONN-2 are displayed in Fig.~\ref{fig:numerical_example_1}. We observe that when the initial shape is circle as in Fig.~\ref{fig:numerical_example_1}(a), both Fireshape and AONN-2 can converge to the reference optimized shape. While, in the case with ellipse shape initialization as in   Fig.~\ref{fig:numerical_example_1}(b), the meshes in the middle of the right side of the computational domain are severely compressed in Fireshape, and a small protrusion occurs eventually. In contrast, AONN-2 stably converges to the reference shape. In the last case, a rectangle that is significantly different from the reference shape is set to the initial shape, which makes the optimization process even more difficult. Due to the restriction of moving mesh, Fireshape fails in this situation. Unlike Fireshape, AONN-2 still roughly converges to the reference shape as shown in Fig.~\ref{fig:numerical_example_1}(c), which indicates that AONN-2 is more flexible with large deformation and more robust to various initial shapes.

\subsection{Pipe optimization constrained by Stokes equations}

\begin{figure}[!htb]
  \centering
  \includegraphics[width=0.36\textwidth]{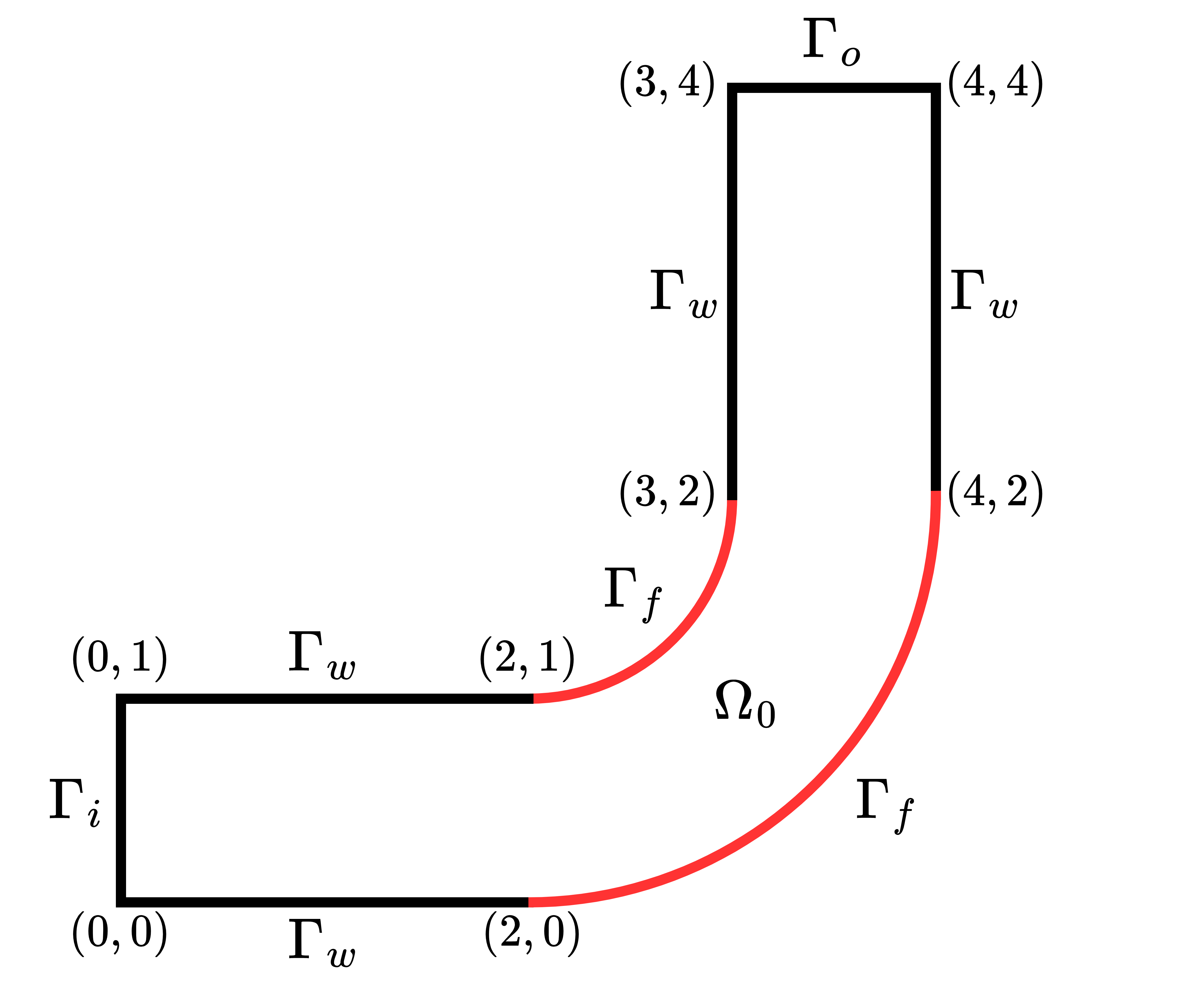}
  \caption{The initial shape of the pipe optimization problem. $\Gamma_{i}$ is inflow boundary and $\Gamma_{o}$ is outflow boundary. $\Gamma_{w}$ and $\Gamma_{f}$ are both no-slip boundaries, and $\Gamma_{f}$ is free boundary which needs to be optimized to minimize the energy dissipation.}
  \label{fig:ELBOW}
\end{figure}

The second problem is for two-dimensional pipe optimization \cite{ schmidt2010efficient, paganini2021fireshape}, which aims to optimize the shape of a pipe to minimize the energy dissipation in the fluid. The objective functional and the constraint due to the Stokes equations are as follows:
\begin{equation}\label{eq:test2}
\left\{
 \begin{aligned}
  &\min\limits_{\mathbf{u}, \Gamma_{f}} J(\mathbf{u},\Omega):=\nu \int_{\Omega}\left\| \nabla \mathbf{u} \right\|_{F}^{2}\,\mathrm{d}\mathbf{x}, \\ 
  &\text{subject to} 
  \left \{
  \begin{aligned}
    -\nabla p + \nu \Delta \mathbf{u}&=\mathbf{0} && \text{ in } \Omega, \\
    \text{div}\,\mathbf{u} &=0 &&  \text{ in } \Omega,\\
    \mathbf{u}&=(u_{in},0)^{\top} && \text{ on } \Gamma_{i}, \\
    \mathbf{u}&=\mathbf{0} &&  \text{ on } \Gamma_{w}\cup \Gamma_{f},\\
    p\mathbf{n}-\nu \partial_\mathbf{n}\mathbf{u}&=\mathbf{0} &&  \text{ on }  \Gamma_{o},\\
    \text{Vol}(\Omega)&=V_{0}.
  \end{aligned}
 \right.
 \end{aligned}
\right.
\end{equation}
where $\nu=1/400$ is the reciprocal of the Reynolds number, and the horizontal velocity profile at the inlet is set to $u_{in}=4(1-x_{2})x_{2}$. The initial domain is illustrated in Fig.~\ref{fig:ELBOW} and the volume of initial domain is set to $V_0$. In this example, only the free boundary $\Gamma_{f}$ can be optimized and the no-slip boundary $\Gamma_{w}$ is  kept fixed. Besides, the volume of the domain needs to remain unchanged during optimization.
Due to the self-adjoint property of problem \eqref{eq:test2}, the adjoint step can be omitted, and we can easily derive the shape derivative:
\begin{equation}
\label{eq:pipe_dJ}
    \mathrm{d}_{\Omega}J(\mathbf{u}, \Omega;\mathbf{V}) = \int_{\Gamma_{f}}\nu
    (\partial_{\mathbf{n}}\mathbf{u})^2\mathbf{n}\cdot \mathbf{V} \mathrm{d}s,
\end{equation}
and regularization equation: 
\begin{equation}
\label{eq:pipe_reg}
\left\{
    \begin{aligned}
      -\Delta \mathbf{\Phi}+\mathbf{\Phi}&=\mathbf{0} && \text{ in } \Omega,\\
      \partial_\mathbf{n}\mathbf{\Phi}+\nu(\partial_{\mathbf{n}}\mathbf{u})^2\mathbf{n}&=\mathbf{0} &&\text{ on }\Gamma_{f},\\
      \mathbf{\Phi}&=\mathbf{0} && \text{ on } \Gamma_{i}\cup \Gamma_{w}\cup \Gamma_{o}.
    \end{aligned}
  \right.
\end{equation}

In order to apply AONN-2 to this problem, two neural networks are used to represent the state variables, namely the velocity field $\mathbf{u}$ and pressure field $p$, and one neural network is used to represent the descent direction $\mathbf{\Phi}$. The network for the velocity field contains 3 residual blocks with 20 neurons per hidden layer, the network for the pressure field contains 3 residual blocks with 15 neurons per layer, and the network for the descent direction contains 3 residual blocks with 20 neurons per hidden layer. These networks take spatial coordinates $(x_1,x_2)$ as the input and  $\mathbf{u},p,\mathbf{\Phi}$ as the output, respectively at the corresponding location. For solving the state equation and regularization equation based on PINNs, $6000$ collocation points inside the domain and 1450 collocation points on the boundary are sampled. Among these boundary collocation points, 450 points on $\Gamma_{f}$ are also used as shape representing points. As a comparison, in Fireshape, the P2-P1 Taylor-Hood finite element method is used for spatial discretization and the Scalable Nonlinear Equations Solvers (SNES) in PETSc are used to solve the discretized system. A total of 1566, 6264 and 25056 triangles are respectively employed to discrete the computational domain.

\begin{table}[!htb]
\caption{The experiment settings and test results for the pipe optimization.}
    \centering
    \begin{small}
    \begin{tabular}{l ccc c}
    \toprule
        & \multicolumn{3}{c}{Fireshape} & \multicolumn{1}{c}{AONN-2}\\
    \cmidrule(lr){2-4} \cmidrule(lr){5-5}
    Initial volume & 6.3562 & 6.3562 & 6.3562 & 6.3562 \\ 
    Optimized volume   &  6.3564 & 6.3556  & 6.3554 & 6.3562 \\
    Initial objective  &  0.0851& 0.0851   &  0.0851  & 0.0851  \\
    Optimized objective & 0.0785  & 0.0785  & 0.0785  & 0.0778 \\
    Shape representation & 1566(triangles)  &  6264(triangles) & 25056(triangles) & 450(points) \\
    Collocation points ($M,N$) & -  &  - & -  &  (1450,6000) \\
    Network parameters & -  & -  & -  & 5665  \\
    Iteration number        &  36 & 48  & 36 & 20 \\
    \bottomrule
    \end{tabular}
    \end{small} %
    \label{tab:test2}
\end{table}

\begin{figure}[!htb]
  \centering
   \subfigure[Fireshape]
  {\includegraphics[width=0.48\textwidth]{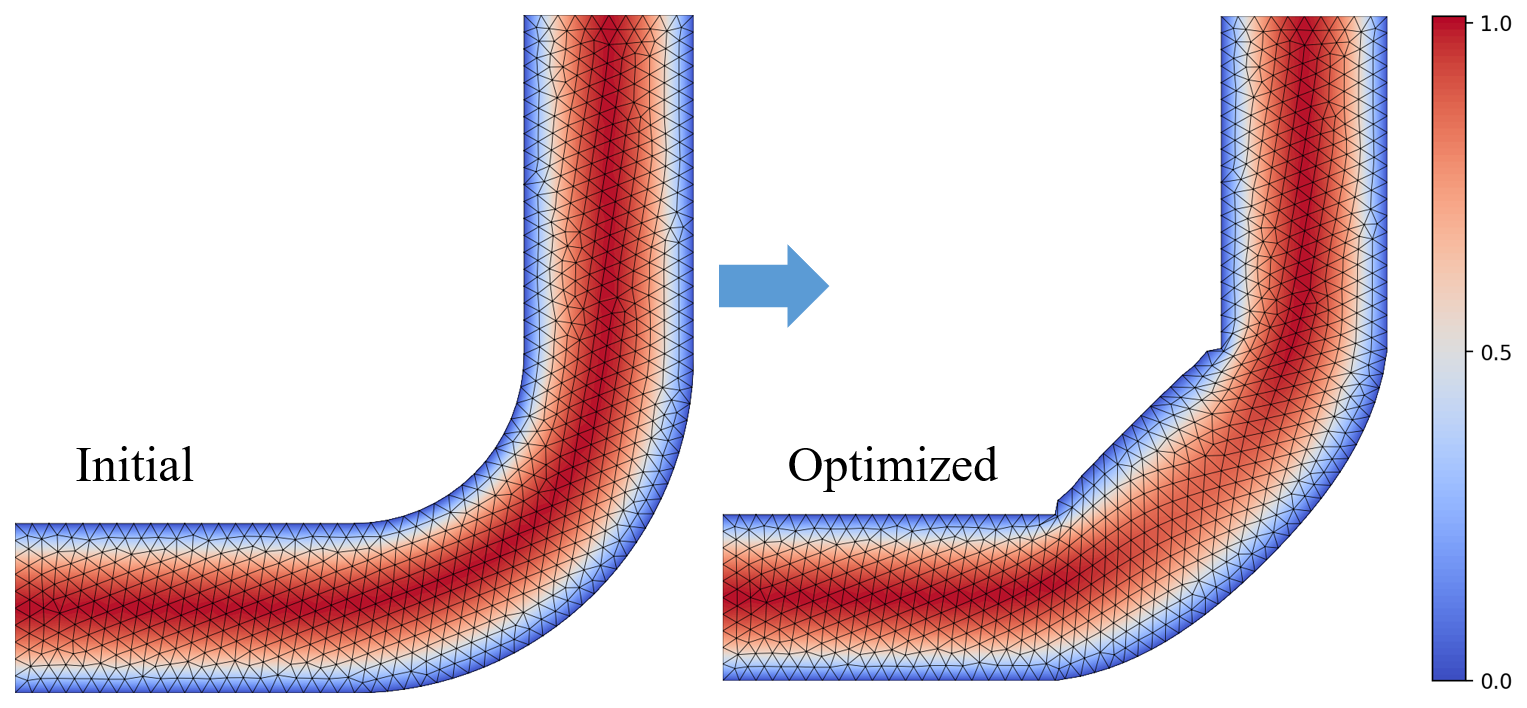}}
  \subfigure[AONN-2]
  {\includegraphics[width=0.48\textwidth]{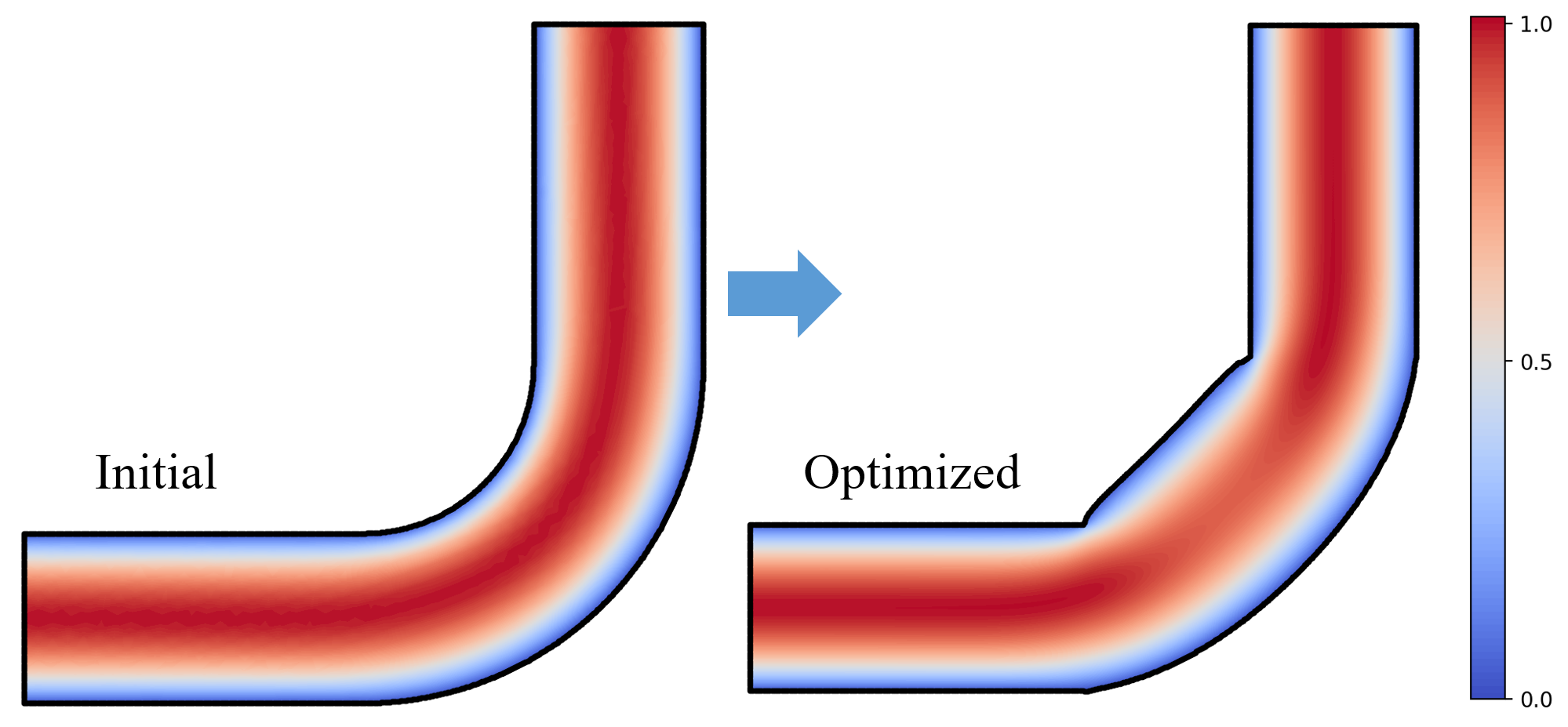}}
  \caption{Pipe optimization results by: (a) Fireshape and (b) AONN-2. The color represents the velocity magnitude.}
  \label{fig:pipe_result}
\end{figure}

The experiment settings and results are listed in Table~\ref{tab:test2}. From the table, we can see that the values of the optimized objective functionals by using Fireshape with different numbers of triangles are almost the same, which are all higher than the value optimized by AONN-2. Besides, AONN-2 keeps the volume of the domain unchanged during the process of optimization, while Fireshape slightly changes the volume of the domain. To make further comparison, we show  the initial and optimized shapes obtained by AONN-2 and Fireshape in Fig.~\ref{fig:pipe_result}.
It can be observed from Fig.~\ref{fig:pipe_result}(a) that, when using Fireshape, the meshes are severely squeezed at the junctions of the free boundary and the fixed boundary, despite of the introduction of regularization terms. And as seen in Fig.~\ref{fig:pipe_result}(b), AONN-2 does not suffer from this problem, thanks to its mesh-free property that can lead to lower optimized value.

\subsection{Obstacle optimization constrained by Stokes equations}
\begin{figure}[!htb]
  \centering
  \subfigure[Case I: one obstacle.]
  {\includegraphics[width=0.49\textwidth]{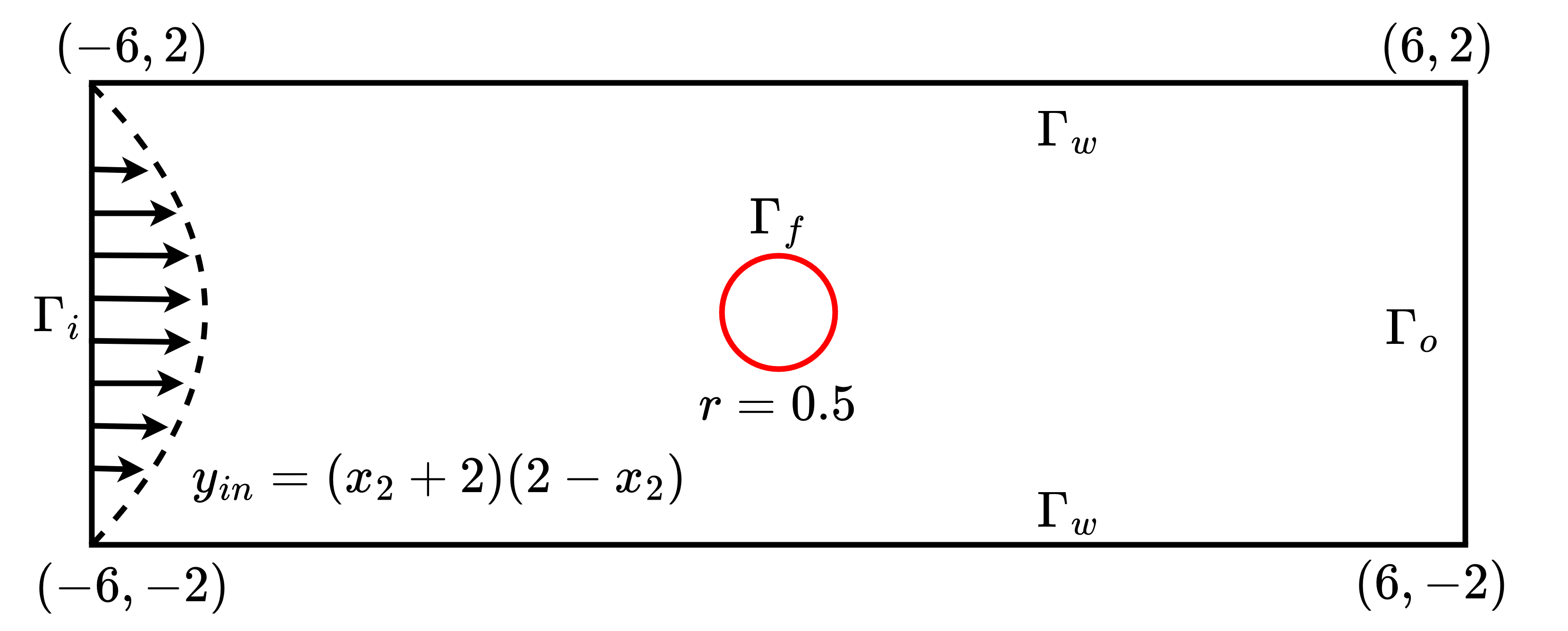}}
  \subfigure[Case II: two obstacles.]
  {\includegraphics[width=0.49\textwidth]{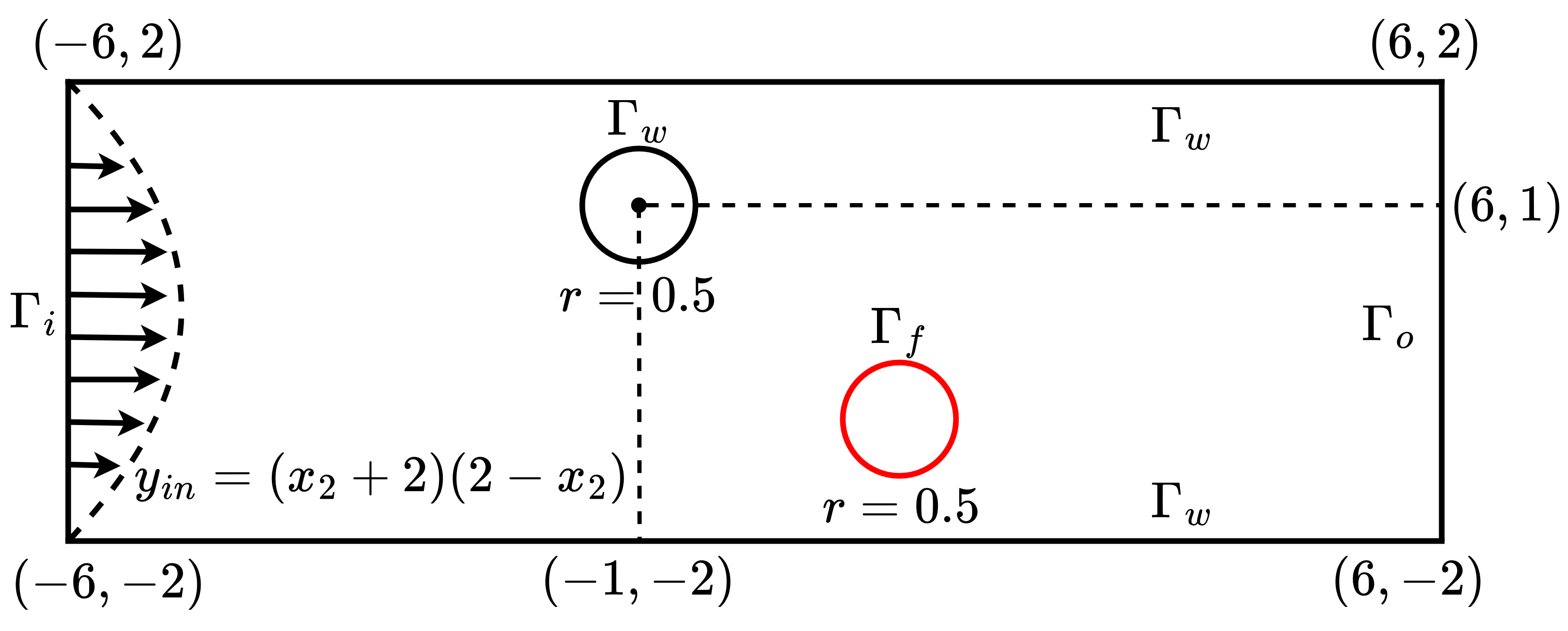}}
  \caption{The initial geometries of two cases in the obstacle optimization problem. $\Gamma_{i}$ is inflow boundary and $\Gamma_{o}$ is outflow boundary. $\Gamma_{w}$ and $\Gamma_{f}$ are no-slip boundaries, where $\Gamma_{f}$ with red color (the red circle) is the free boundary that needs to be optimized.}
\label{fig:obstacle}
\end{figure}

Next, we consider the obstacle optimization constrained by Stokes equations in a pipe flow, which can be frequently found in aerodynamic applications. The aim of the obstacle optimization is to search for the shape that minimizes the drag on the obstacle, which is equivalent to search for the shape that minimizes the energy dissipation in the fluid due to the shear forces \cite{paganini2018higher, bello1997differentiability}. 
We take the same objective functional and PDE constraint as the pipe optimization given in equation \eqref{eq:test2}, with a viscosity value of $\nu=1/80$. This setting leads to the same formulations of shape derivative and regularization equation.
Two different scenarios of the obstacle are established as in Fig.~\ref{fig:obstacle}. In the first case, the initial shape of the obstacle is a circle with a radius of $0.5$, placed in the center of the pipe. In the second case, the initial obstacle is located below the right of the pipe center, and there is a fixed circular obstacle at the symmetrical position of the center. In these two cases, the barycenter of the obstacle is not fixed. The geometry parameters of the computational domain are marked in Fig.~\ref{fig:obstacle}.

To implement the AONN-2 algorithm, three neural networks are employed to represent the velocity field $\mathbf{u}$, the pressure field $p$, and the descent direction $\mathbf{\Phi}$, respectively, which are all comprised of $2$ residual blocks with $15$ neurons per hidden layer. Different numbers of collocation points and shape representing points are adopted, which are listed in Table~\ref{tab:test3}. For comparison, Fireshape is also used to solve the problem, with the P2-P1 Taylor-Hood finite element method as the spatical discretization and the Krylov subspace method in PETSc as the linear solver. We use the number of triangles as the indicator of shape representation of finite element. The corresponding experiment settings and results are listed in Table~\ref{tab:test3}.

\begin{table}[!htb]
\caption{The experiment settings and test results for the obstacle optimization.}
    \centering
    \begin{scriptsize}
    \begin{tabular}{l cccc cccc}
    \toprule
      Case I  & \multicolumn{4}{c}{Fireshape} & \multicolumn{4}{c}{AONN-2}\\
    \cmidrule(lr){2-5} \cmidrule(lr){6-9}
    Initial objective  & $0.4574$ &  $0.4574$  & $0.4574$   & $0.4574$ & $0.4574$&$0.4574$& $0.4574$& $0.4574$ \\
    Optimized objective  & $0.4264$  & $0.4260$  & $0.4259$  & $0.4258$ & $0.4187$ & $0.4190$ & $0.4188$ & $0.4190$ \\
    Shape representation & 872  &  3488 & 13952 & 55808 & 600 & 600 & 1200 & 1200 \\
    Collocation points ($M,N$) & -  &  - & -  &  - & (3800,6000) & (3800,6000)& (4400,12000) & (4400,12000)\\
    Network parameters & -  & -  & -  & - & 3172 & 5092 & 3172 & 5092\\
    Iteration number        &  70 & 83  & 68 & 70 & 30 & 30& 30& 30\\
    \cmidrule(lr){1-9}
      Case II & \multicolumn{4}{c}{Fireshape} & \multicolumn{4}{c}{AONN-2}\\
    \cmidrule(lr){2-5} \cmidrule(lr){6-9}
    Initial objective  & $0.5507$ &  $0.5507$  & $0.5507$   & $0.5507$ & $0.5507$&$0.5507$& $0.5507$& $0.5507$ \\
    Optimized objective   & $0.4751$  & $0.4744$  & $0.4742$  & $0.4744$ & $0.4121$ & $0.3676$ & $0.3810$ & $0.3638$ \\
    Shape representation & 1404  &  5616 & 22464 & 89856 & 600 & 600 & 1200 & 1200 \\
    Collocation points ($M,N$) & -  &  - & -  &  - & (4400,6000) & (4400,6000)& (5000,12000) & (5000,12000)\\
    Network parameters & -  & -  & -  & - & 3172 & 5092 & 3172 & 5092\\
    Iteration number        &  68 & 78  & 86 & 88 & 100 & 100& 100& 100\\
    \bottomrule
    \end{tabular}
    \end{scriptsize} %
    \label{tab:test3}
\end{table}

\begin{figure}[!htb]
  \centering
  \subfigure[I: initial (Fireshape)]
  {\label{fig:obstacle-flow-1a}\includegraphics[width=0.46\textwidth]{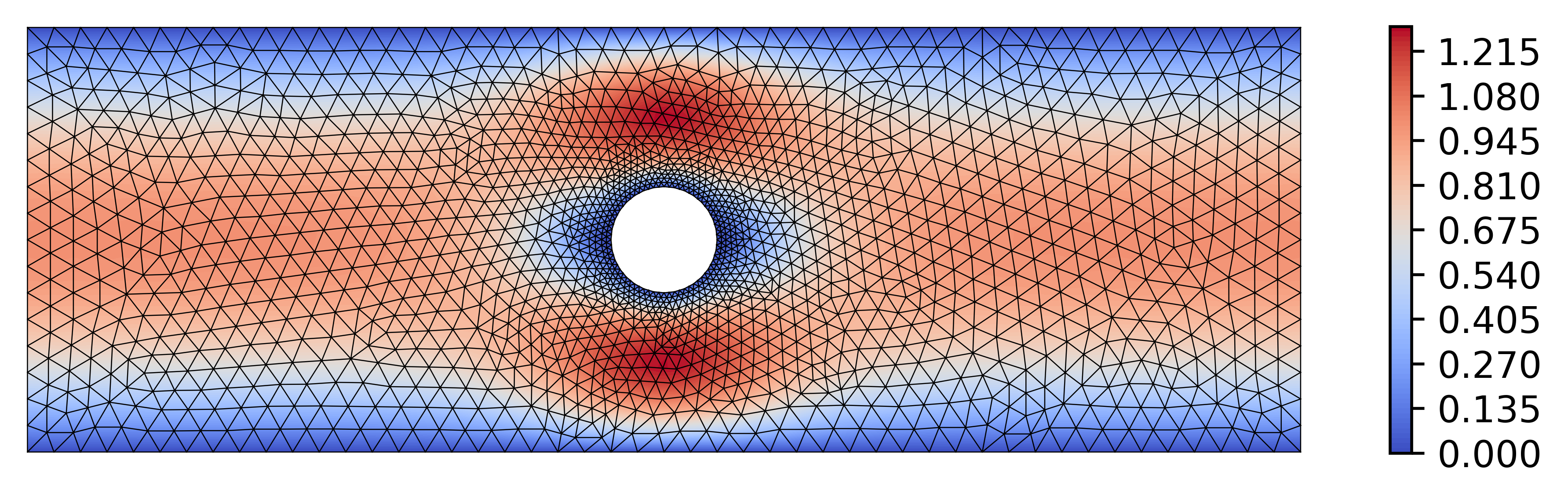}}
  \subfigure[I: initial (AONN-2)]
  {\label{fig:obstacle-flow-1b}\includegraphics[width=0.46\textwidth]{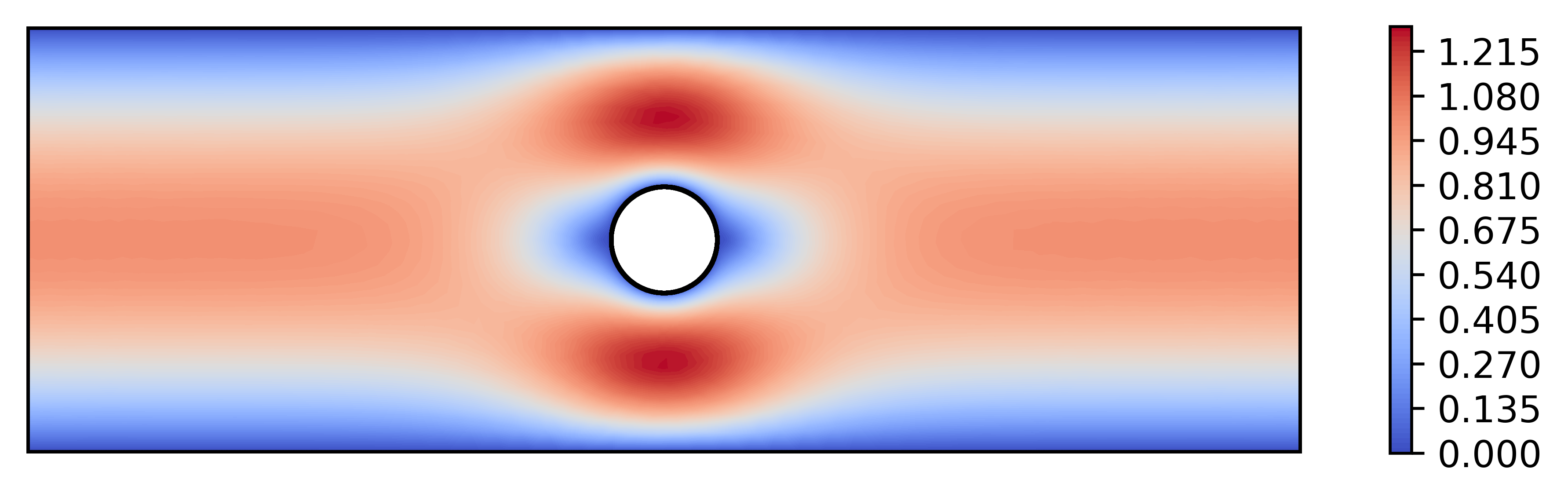}}
  \subfigure[I: optimized (Fireshape)]
  {\label{fig:obstacle-flow-1c}\includegraphics[width=0.46\textwidth]{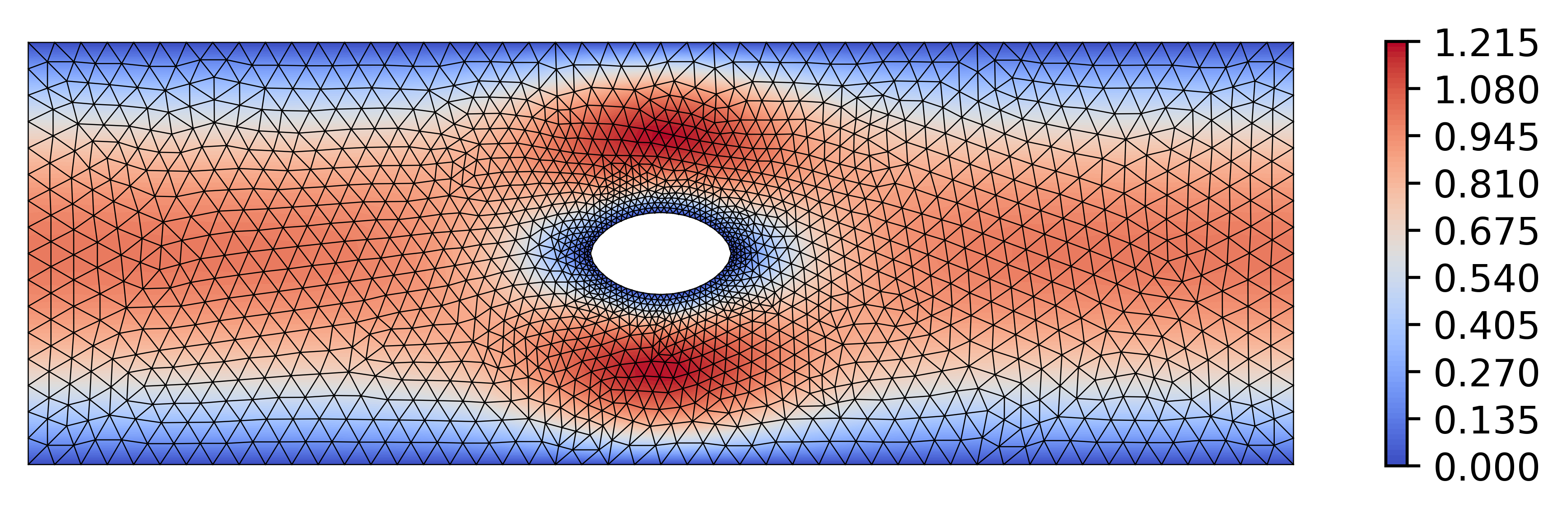}}
  \subfigure[I: optimized (AONN-2)]
  {\label{fig:obstacle-flow-1d}\includegraphics[width=0.46\textwidth]{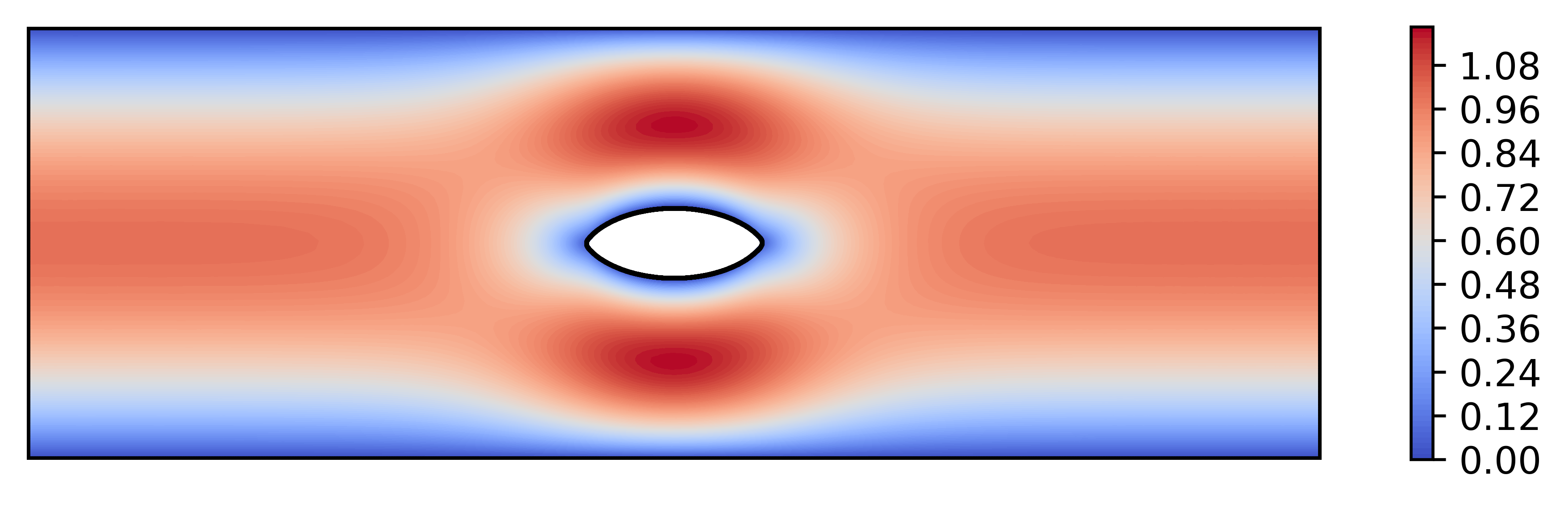}}
  \caption{Obstacle optimization results for test case I by Fireshape and AONN-2. The color represents the velocity magnitude.}
\label{fig:obstacle-flow-1}
\end{figure}

\begin{figure}[!htb]
  \centering
  \subfigure[II: initial (Fireshape)]
  {\label{fig:obstacle-flow-2a}\includegraphics[width=0.46\textwidth]{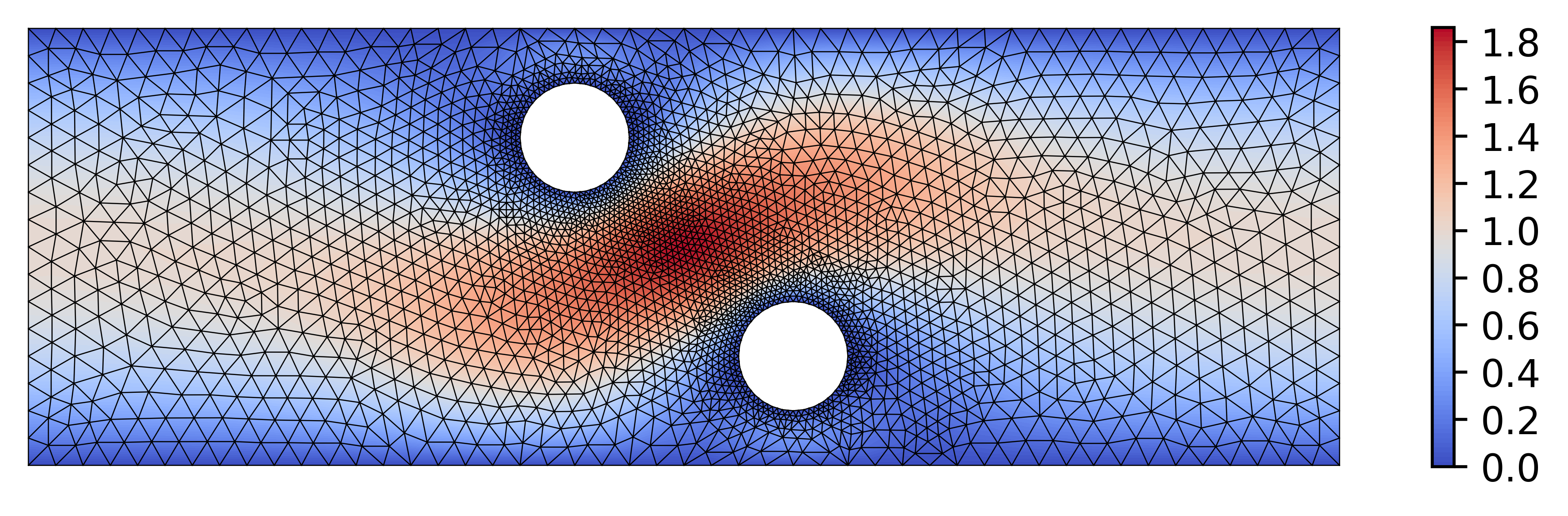}}
  \subfigure[II: initial (AONN-2)]
  {\label{fig:obstacle-flow-2b}\includegraphics[width=0.46\textwidth]{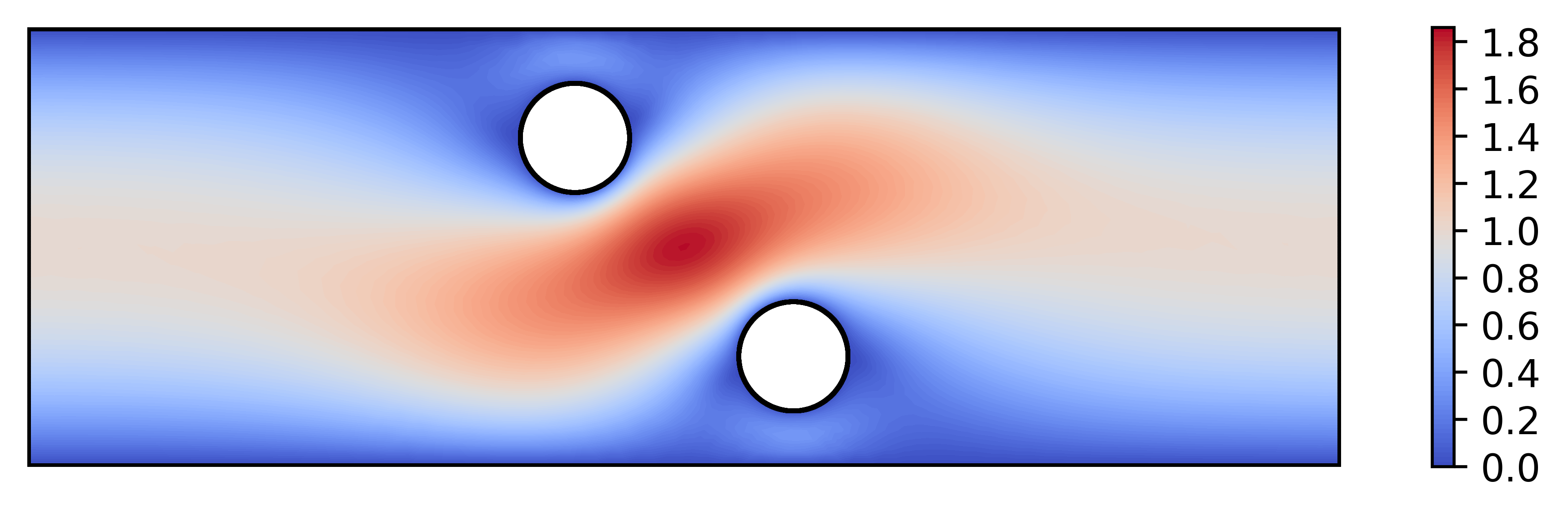}}
  \subfigure[II: optimized (Fireshape)]
  {\label{fig:obstacle-flow-2c}\includegraphics[width=0.46\textwidth]{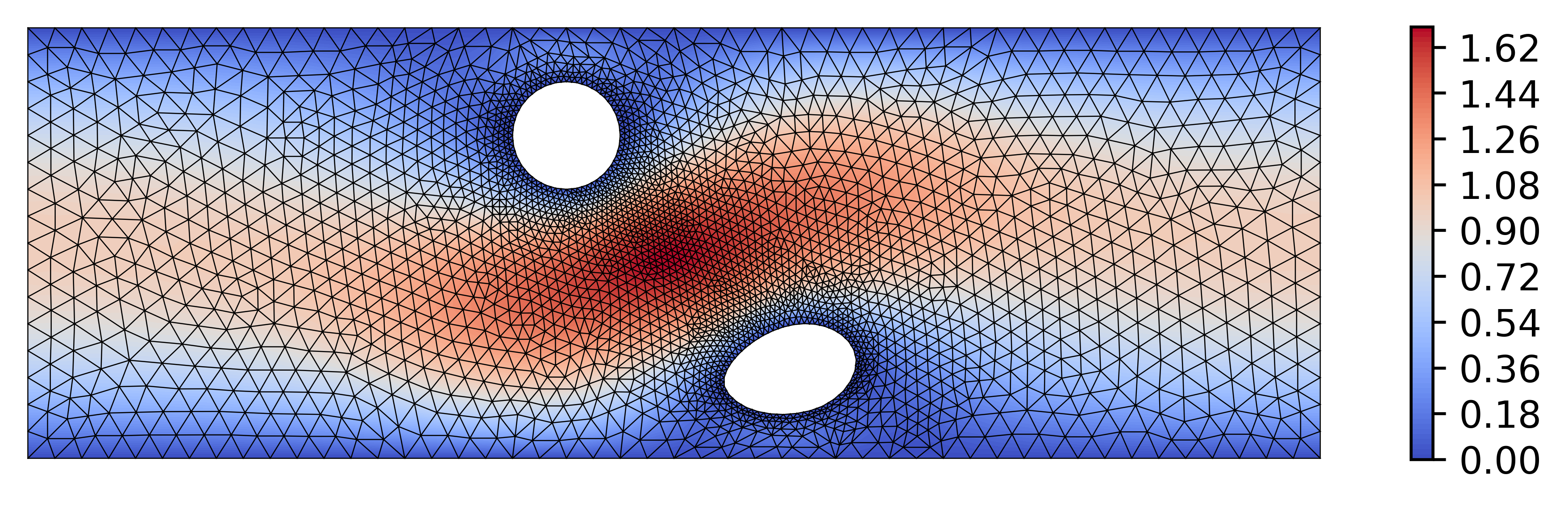}}
  \subfigure[II: optimized (AONN-2)]
  {\label{fig:obstacle-flow-2d}\includegraphics[width=0.46\textwidth]{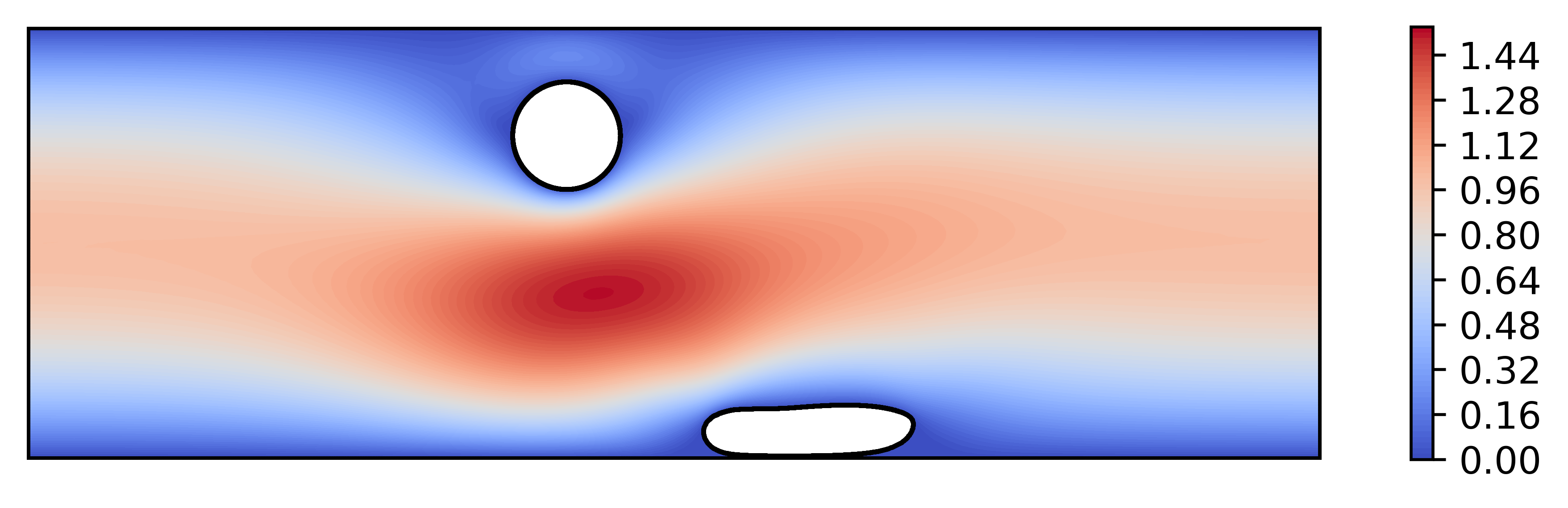}}
  \caption{Obstacle optimization results for test case II by Fireshape and AONN-2. The color represents the velocity magnitude.}
\label{fig:obstacle-flow-2}
\end{figure}

In Fig.~\ref{fig:obstacle-flow-1} and \ref{fig:obstacle-flow-2} we show the flow fields before and after the optimization, for test case I and II, respectively. By comparing Fig.~\ref{fig:obstacle-flow-1c} and \subref{fig:obstacle-flow-1d}, we can observe that AONN-2 can outperform Fireshape and obtain an optimized obstacle with flatter shape, leading to lower value of the objective functional. And from Fig.~\ref{fig:obstacle-flow-2c} and \subref{fig:obstacle-flow-2d}, it can be found that the position of optimal shape of the obstacle obtained by AONN-2 is relatively close to the boundary, which makes the flow smoother and thus can reduce the energy dissipation. While in Fireshape, the quality of surrounding meshes of the obstacle needs to be guaranteed (e.g. adding regularization terms), which limits the deformation of the obstacle itself. In contrast, the shape deformation in AONN-2 is more flexible, and thus lower value of the objective functional is attained.

\subsection{Channel optimization constrained by Naiver-Stokes equations}

\begin{figure}[!htb]
  \centering
  \vspace{-12pt}
  \subfigcapskip=-2pt
  \subfigure[Case I: convex initialization]
  {\includegraphics[width=0.33\textwidth]{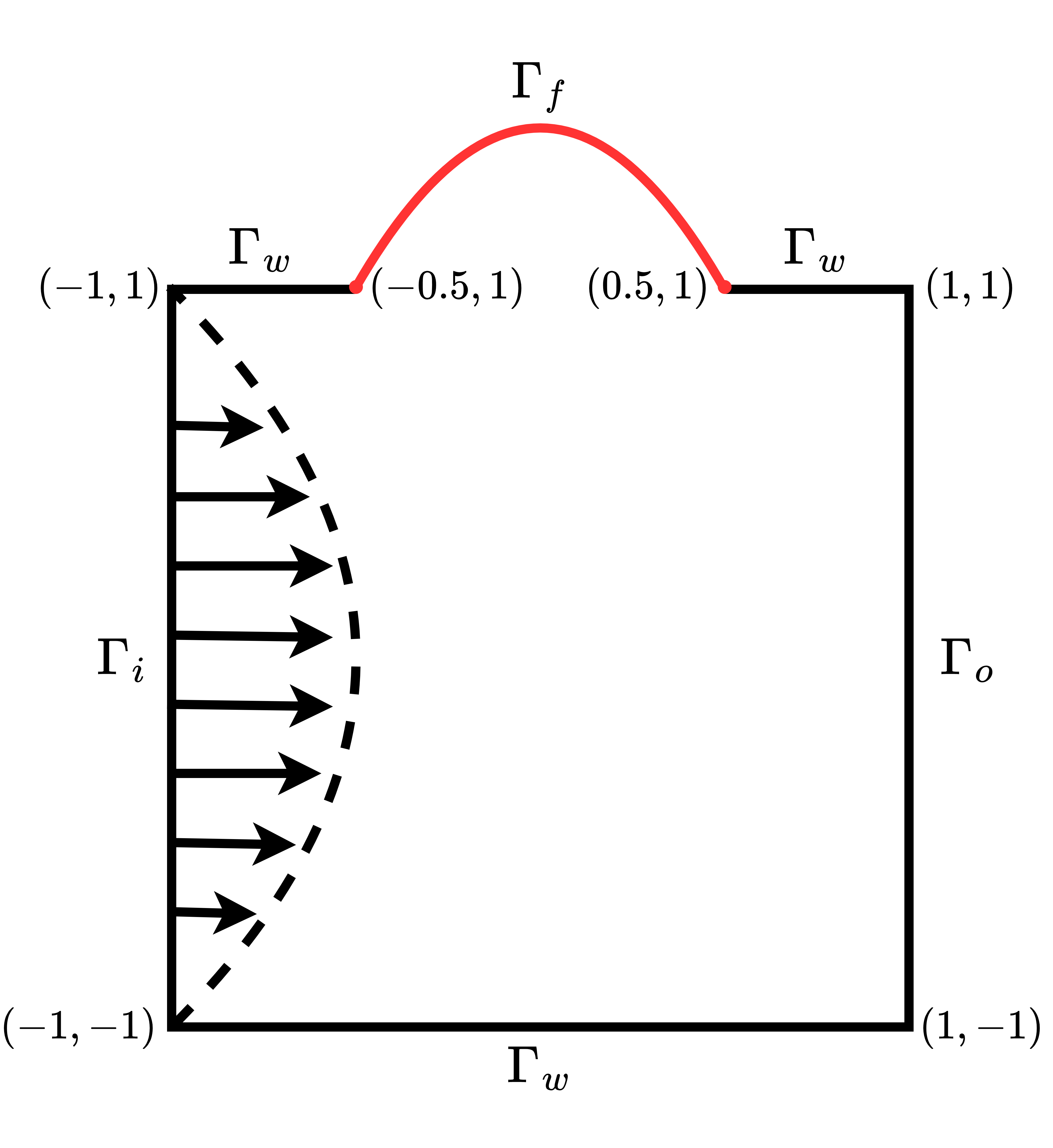}}
  \hspace{20pt}
  \subfigure[Case II: concave initialization]
  {\includegraphics[width=0.33\textwidth]{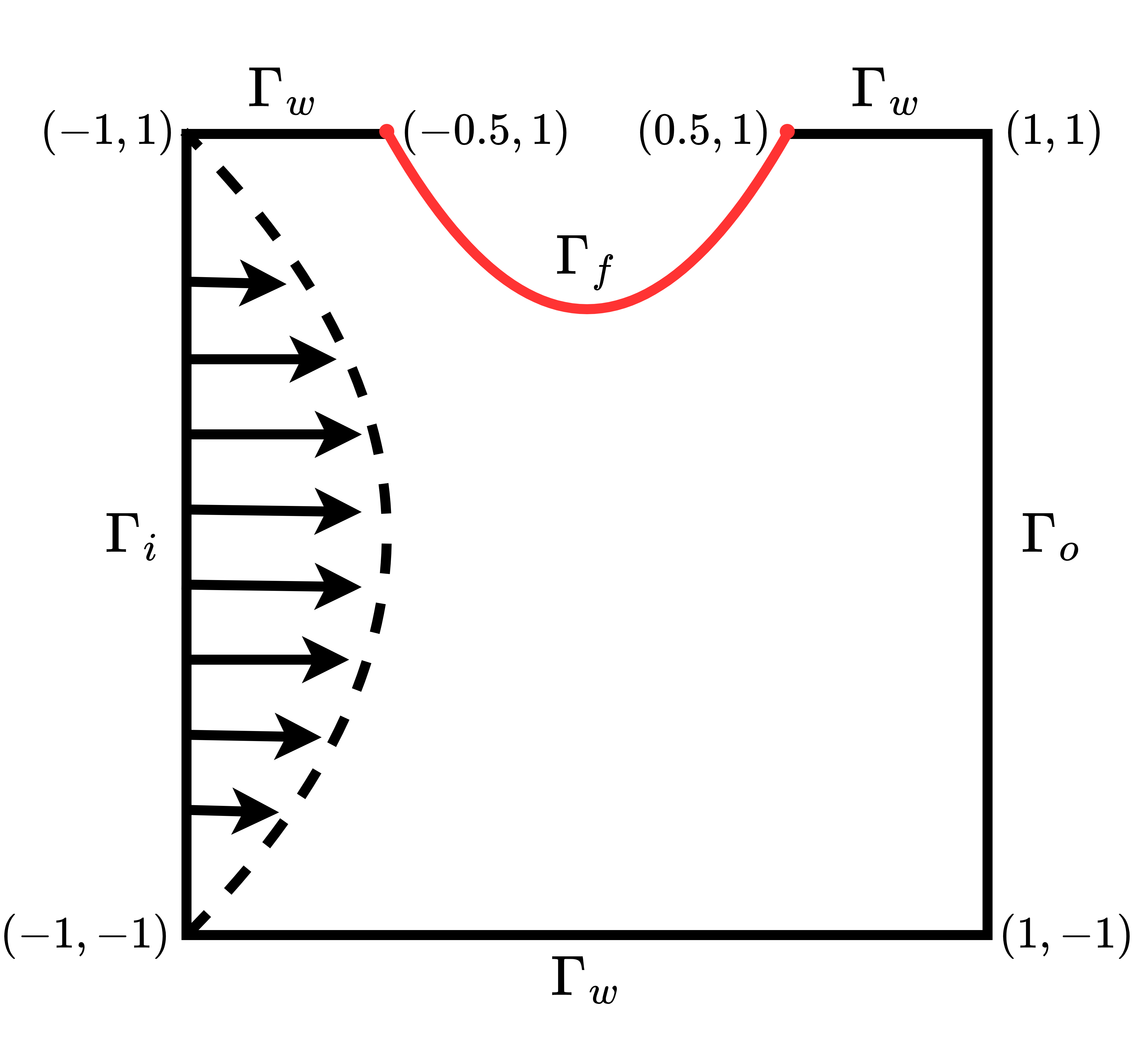}}
  \caption{Two initial shapes in the channel optimization problem. A prescribed horizontal velocity profile $\mathbf{u}=(u_{in},0)^{\top}$ is assigned on the inlet boundary $\Gamma_{i}$. The flow leaves the domain at the outflow boundary $\Gamma_{o}$, and the remaining boundaries $\Gamma_{w}$, $\Gamma_{f}$ are no-slip walls, where $\Gamma_{f}$ needs to be optimized to minimize the objective functional.}
  \label{fig:BUMP_2D}
\end{figure}

As the final example, the shape optimization constrained by the 2D incompressible Navier-Stokes equations \cite{kasumba2012vortex} is studied. In this example, an L2-tracking type objective functional is minimized as follows:

\begin{equation}
   \label{eq:bump}
 \left\{
   \begin{aligned}
  &\min\limits_{\mathbf{u}, \Gamma_{f}} J(\mathbf{u},\Omega):=\int_{\Omega}\left\| \mathbf{u} -\mathbf{u}_{d}\right\|^{2}\,\mathrm{d}\mathbf{x} ,\\ 
  &\text{subject to} 
  \left \{
  \begin{aligned}
    -\nabla p - (\mathbf{u}\cdot \nabla)\mathbf{u} + \nu \Delta \mathbf{u}&=\mathbf{0} && \text{ in } \Omega, \\
    \text{div}\,\mathbf{u} &=0 &&  \text{ in } \Omega,\\
    \mathbf{u}&=(u_{in},0)^{\top} && \text{ on } \Gamma_{i}, \\
    \mathbf{u}&=\mathbf{0} &&  \text{ on }  \Gamma_{w}\cup \Gamma_{f},\\
    p\mathbf{n}-\nu \partial_\mathbf{n}\mathbf{u}&=\mathbf{0} &&  \text{ on }  \Gamma_{o},
  \end{aligned}
 \right.
 \end{aligned}
\right.
\end{equation}
where $\nu=1/50$ is the reciprocal of the Reynolds number, $u_{in}(x_{1},x_{2})=2.5(1+x_{2})(1-x_{2})$, and $\mathbf{u}_{d}=(u_{in},0)^{\top}$. Two initial shapes are constructed based on $\overline{\Omega}=[-1,1]\times[-1,1]$ with one piece of its boundary $\{(x_{1},x_{2})|-0.5\leq x_{1}\leq 0.5,x_{2}=1\}$ being replaced by two Bezier curves as shown in Fig.\ref{fig:BUMP_2D}(a) and (b), and the optimal shape corresponding to both initial shapes are exactly the same square domain $\overline{\Omega}$. Under such situation, this problem can be used to test whether the shape optimization method has the ability to turn the curve on the top to a straight line. 

According to problem \eqref{eq:bump}, the adjoint variables $\boldsymbol{\lambda}$ and ${q}$ are defined by the following adjoint equation:

\begin{equation}
  \label{eq:NS adjoint}
    \left\{
    \begin{aligned}
    -\nu\Delta \boldsymbol{\lambda}-\nabla\boldsymbol{\lambda}\cdot \mathbf{u}+(\nabla\mathbf{u})^{\top}\boldsymbol{\lambda}+\nabla q &= \mathbf{u}-\mathbf{u}_{d} && \text{ in } \Omega, \\ 
      \nabla\cdot \boldsymbol{\lambda}&=0 && \text{ in } \Omega,\\
      q\mathbf{n}-\nu \partial_{\mathbf{n}}\boldsymbol{\lambda}-(\mathbf{u}\cdot\mathbf{n})\boldsymbol{\lambda}&=\mathbf{0} && \text{ on } \Gamma_{o},\\
      \boldsymbol{\lambda} &= \mathbf{0} && \text{ on } \Gamma_{w}\cup\Gamma_{f}\cup\Gamma_{i}. \\
    \end{aligned}
  \right.
\end{equation}
Once the state and adjoint variables are obtained, the descent direction $\mathbf{\Phi}$ can be calculated by solving the following regularization equation:

\begin{equation}
  \left\{
    \begin{aligned}
      -\Delta \mathbf{\Phi} + \mathbf{\Phi} &= \mathbf{0}&& \text{ in } \Omega,\\
      \partial_{\mathbf{n}}\mathbf{\Phi} + \nabla \hat{J} \mathbf{n}&= \mathbf{0}  && \text{ on } \Gamma_{f}, \\ 
      \mathbf{\Phi} &= \mathbf{0} && \text{ on }  \Gamma_{i}\cup \Gamma_{w}\cup \Gamma_{o},
    \end{aligned}
  \right.
\end{equation}
where
\begin{equation}
\nabla \hat{J}=\frac{1}{2}\|\mathbf{u}-\mathbf{u}_{d}\|^{2}+\nu (\partial_{\mathbf{n}}\mathbf{u})\cdot(\partial_{\mathbf{n}}\boldsymbol{\lambda}).
\end{equation}

In AONN-2, three neural networks that contain 3 residual blocks with 20 neurons per hidden layer are employed to represent the velocity field $\mathbf{u}$, the adjoint variable $\boldsymbol{\lambda}$ and the descent direction $\mathbf{\Phi}$, respectively, and two neural networks that contain 3 residual blocks with 15 neurons per hidden layer are employed to represent the pressure field $p$ and the adjoint variable $q$, respectively. For training the neural networks to approximate the solutions of the state, adjoint  and regularization equations, $424$ and $849$ points on the domain boundary, and $4000$ and $8000$ points inside the computational domain are sampled for the two cases, respectively. Among these boundary collocation points, $75$ and $150$ points on $\Gamma_{f}$ are used as the shape representing points for each of the two cases. When using Fireshape, the P2-P1 Taylor-Hood finite element method is applied for spatial discretization and the resultant nonlinear equations are solved with PETSc. The corresponding experiment settings and results are listed in Table~\ref{tab:test4}.


\begin{table}[!htb]
 \caption{The experiment settings and results for the channel optimization.}
    \centering
    \begin{scriptsize}
    \begin{tabular}{l cccc cccc}
    \toprule
      Case I  & \multicolumn{4}{c}{Fireshape} & \multicolumn{4}{c}{AONN-2}\\
    \cmidrule(lr){2-5} \cmidrule(lr){6-9}
    Initial objective  & $6.78$ &  $6.78$  & $6.78$   & $6.78$ & $6.78$&$6.78$& $6.78$& $6.78$ \\
    Optimized objective  & $0.0036$  & $0.0031$  & $0.0025$  & $0.0024$ & $0.0017$ & $0.00068$ & $0.0010$ & $0.0011$ \\
    Shape representation & 256  &  1024 & 16384 & 65536 & 75 & 75 & 150 & 150 \\
    Collocation points ($M,N$) & -  &  - & -  &  - & (424,4000) & (424,4000)& (849,8000) & (849,8000)\\
    Network parameters & -  & -  & -  & - & 5013 & 9128 & 5013 & 9128\\
    Iteration number        &  4 & 4  & 5 & 5 & 10 & 10& 10& 10\\
    \cmidrule(lr){1-9}
      Case II & \multicolumn{4}{c}{Fireshape} & \multicolumn{4}{c}{AONN-2}\\
    \cmidrule(lr){2-5} \cmidrule(lr){6-9}
    Initial objective & $6.60$ &  $6.60$  & $6.60$   & $6.60$ & $6.60$&$6.60$& $6.60$& $6.60$ \\
    Optimized objective  & $0.0023$  & $0.0021$  & $0.0018$  & $0.0019$ & $0.00033$ & $0.00031$ & $0.00032$ & $0.00026$ \\
    Shape representation & 288  &  1152 & 4608 & 18432 & 75 & 75 & 150 & 150 \\
    Collocation points ($M,N$) & -  &  - & -  &  - & (424,4000) & (424,4000)& (849,8000) & (849,8000)\\
    Network parameters & -  & -  & -  & - & 5013 & 9128 & 5013 & 9128\\
    Iteration number        &  12 & 12  & 16 & 14 & 20 & 20& 20& 20\\
    \bottomrule
    \end{tabular}
    \end{scriptsize} %
   
    \label{tab:test4}
\end{table}

\begin{figure}[!htb]
  \centering
  \subfigure[\scriptsize{I: initial  (Fireshape)}]
  {\includegraphics[width=0.24\textwidth]{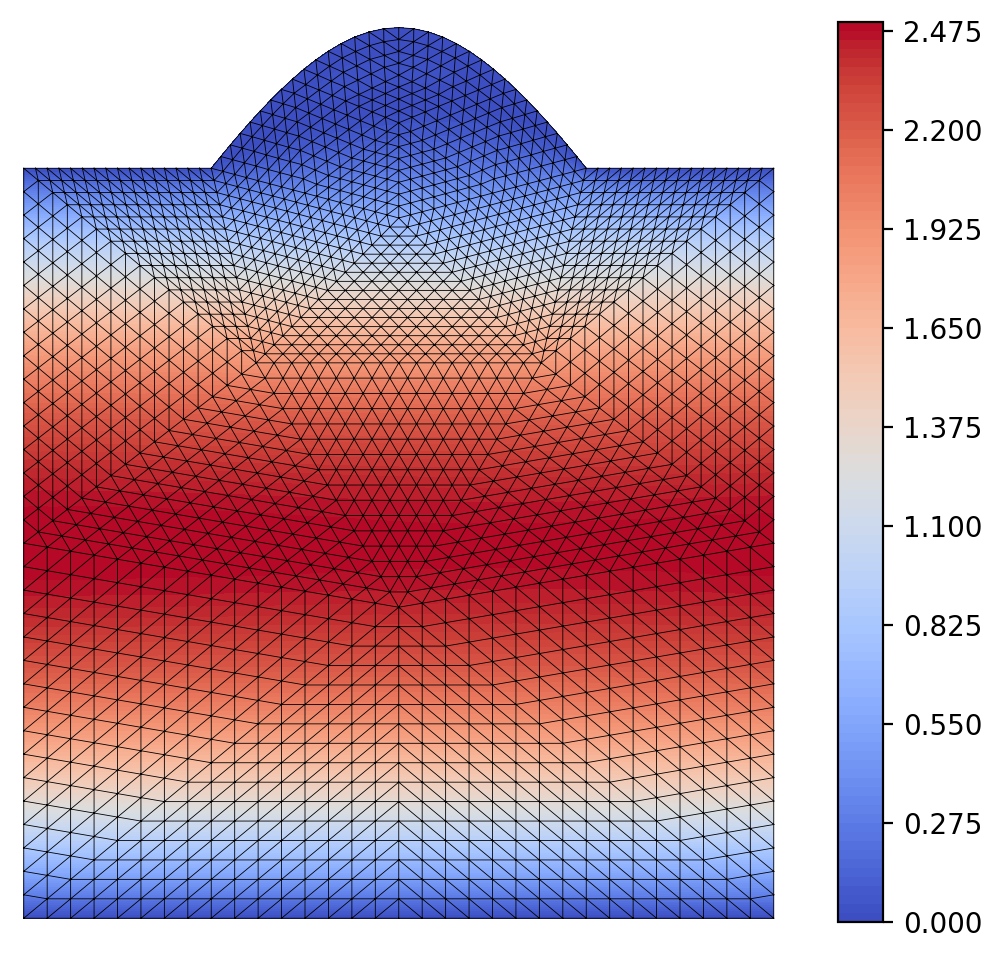}}
   \subfigure[\scriptsize{I: optimized (Fireshape)}]
  {\includegraphics[width=0.24\textwidth]{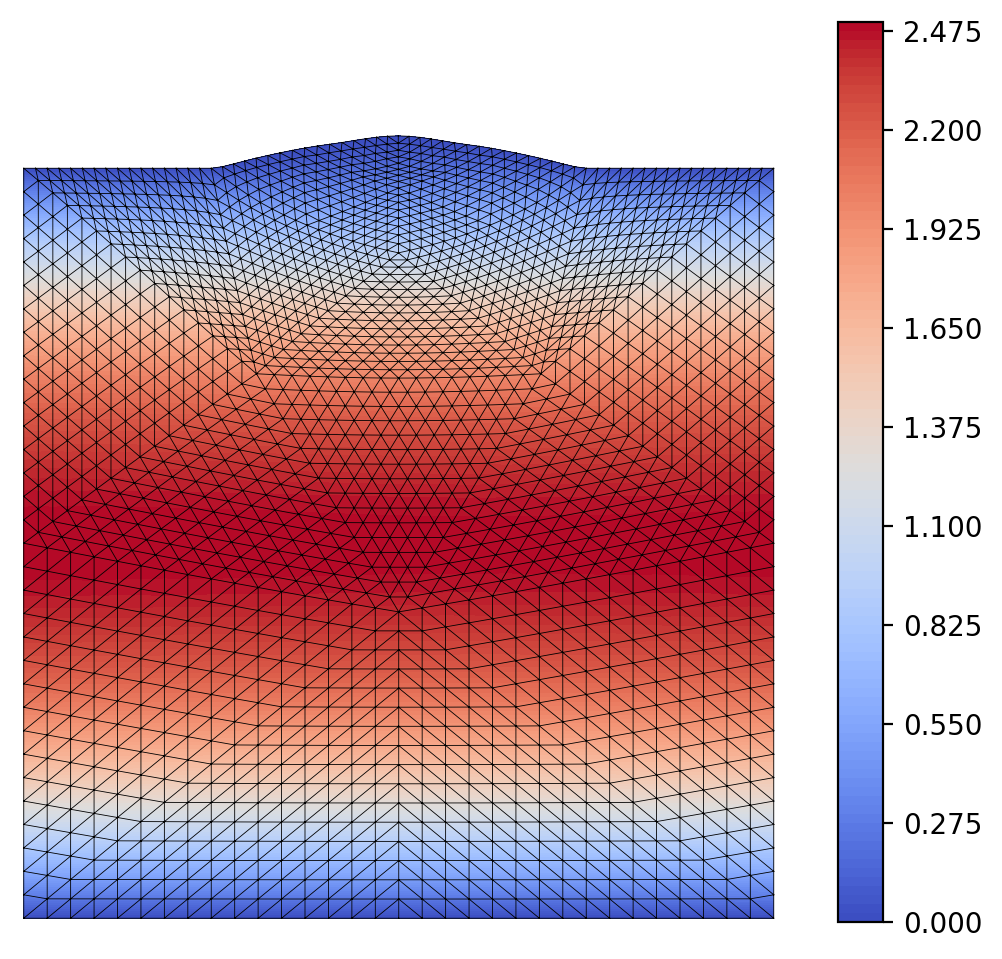}}
  \subfigure[\scriptsize{I: initial (AONN-2)}]
  {\includegraphics[width=0.24\textwidth]{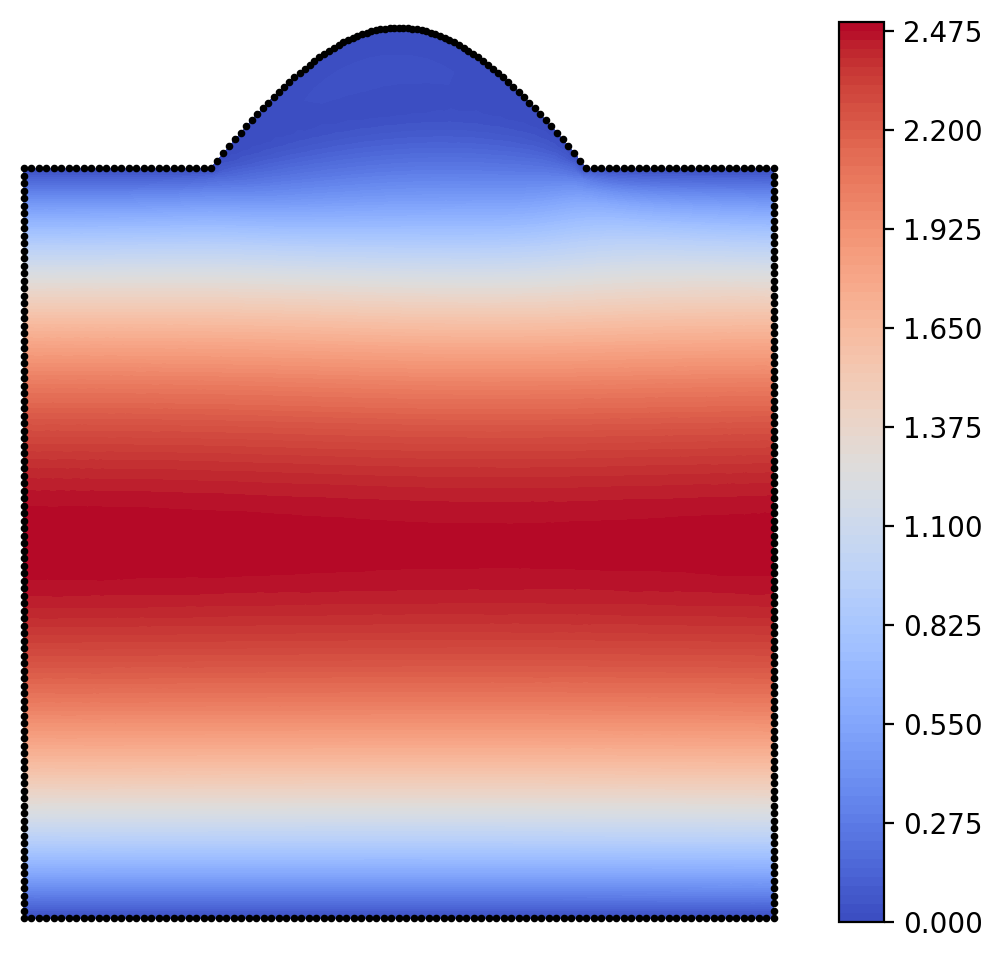}}
  \subfigure[\scriptsize{I: optimized (AONN-2)}]
  {\includegraphics[width=0.24\textwidth]{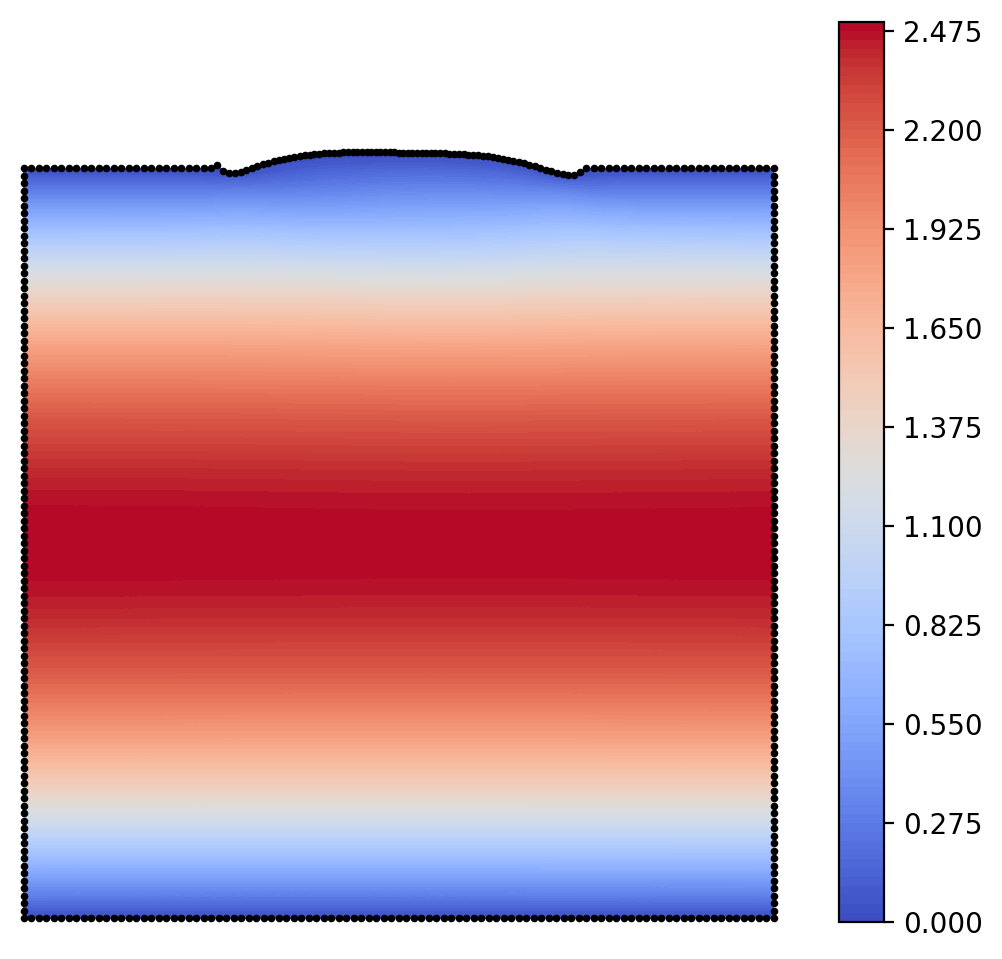}}
  \subfigure[\scriptsize{II: initial (Fireshape)}]
  {\includegraphics[width=0.24\textwidth]{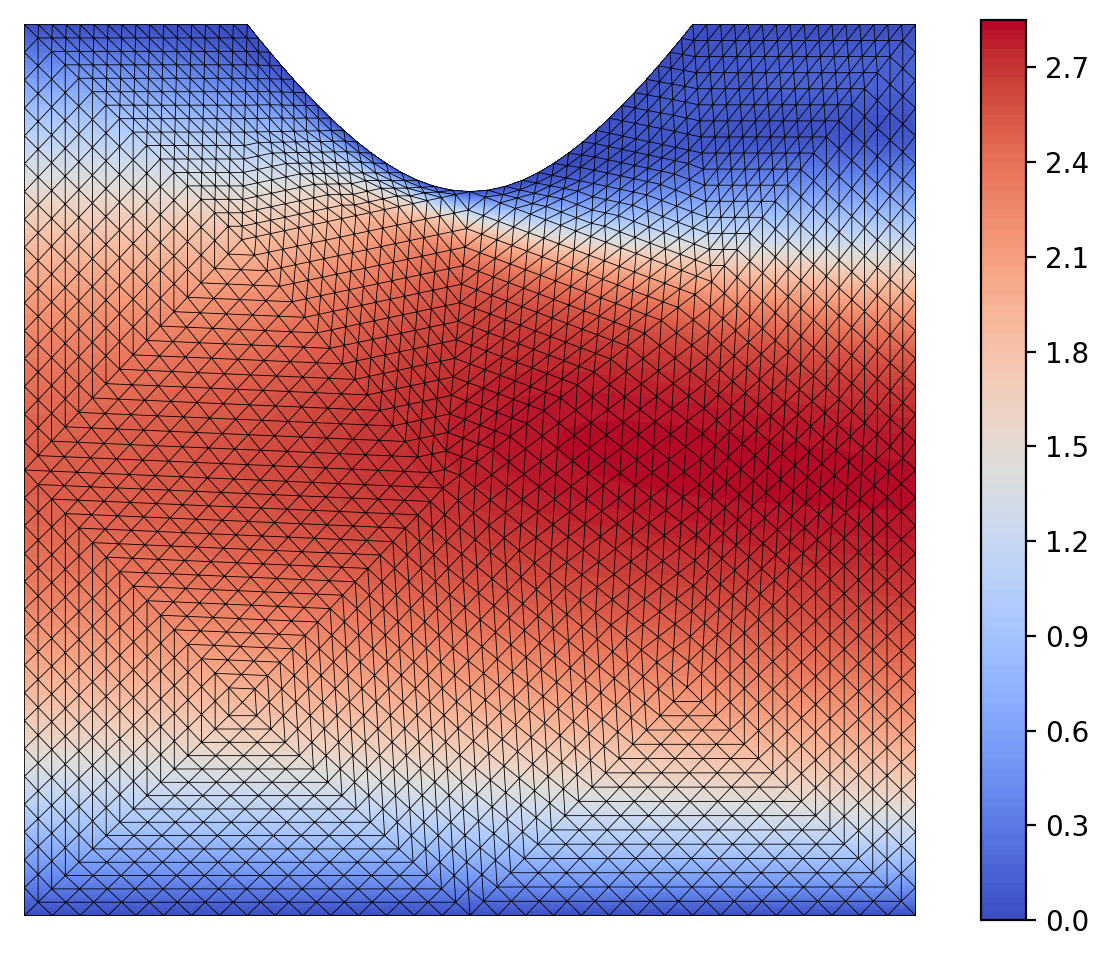}}
   \subfigure[\scriptsize{II: optimized (Fireshape)}]
  {\includegraphics[width=0.24\textwidth]{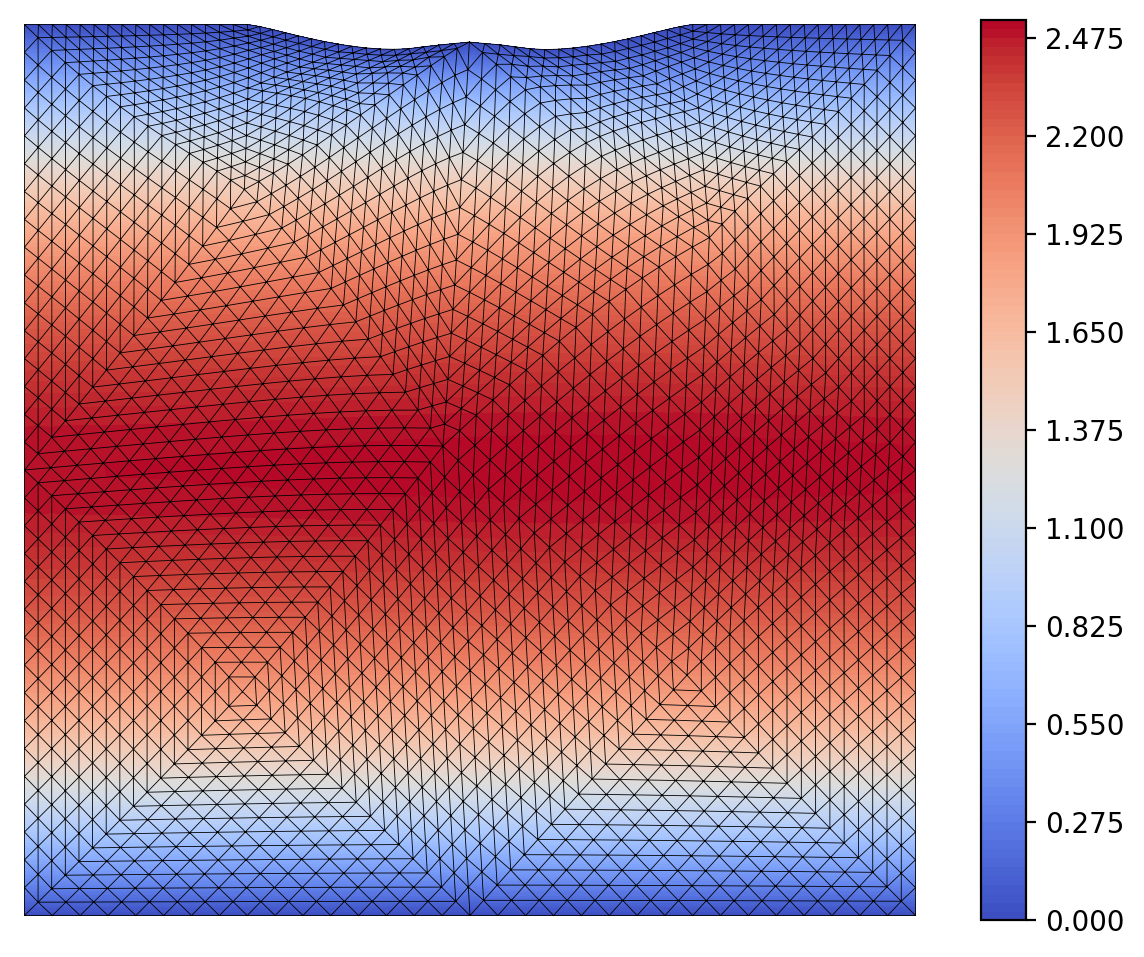}}
  \subfigure[\scriptsize{II: initial (AONN-2)}]
  {\includegraphics[width=0.24\textwidth]{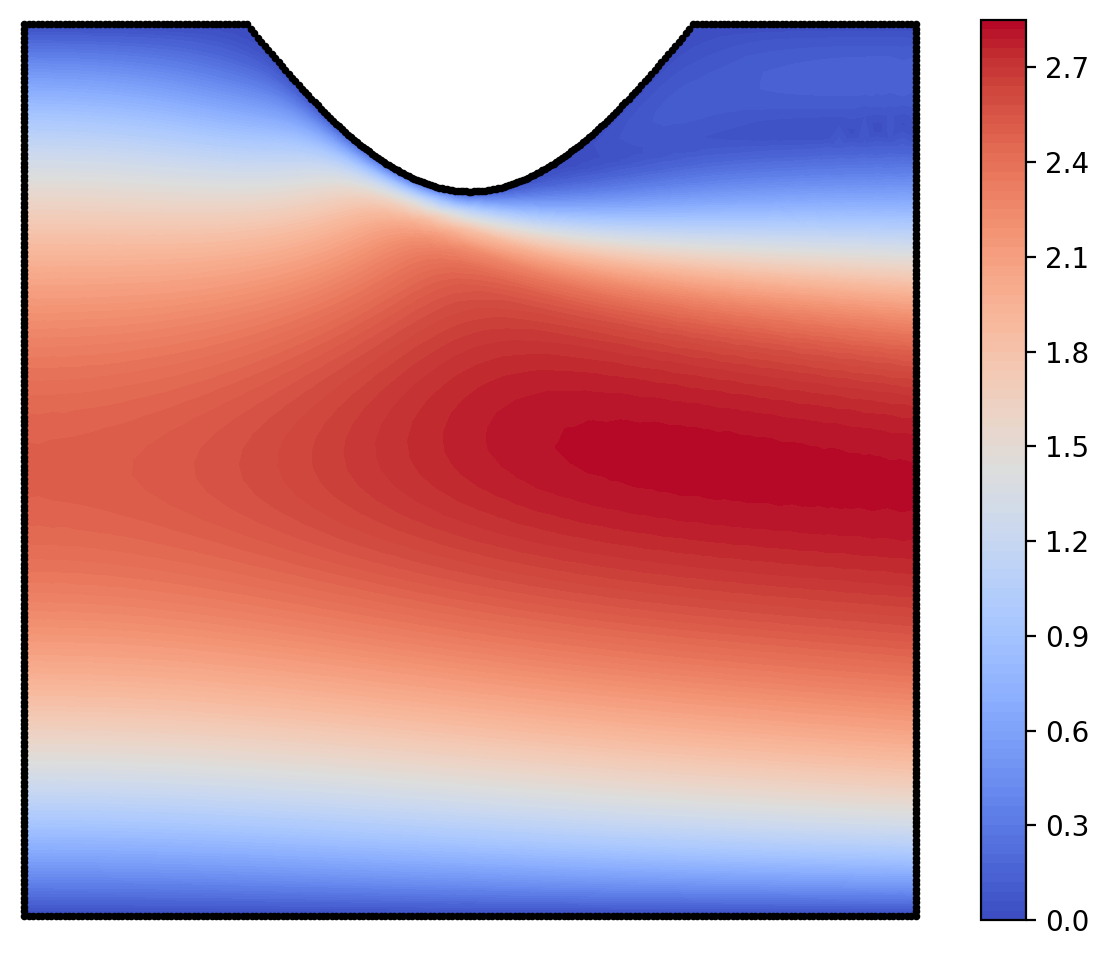}}
  \subfigure[\scriptsize{II: optimized (AONN-2)}]
  {\includegraphics[width=0.24\textwidth]{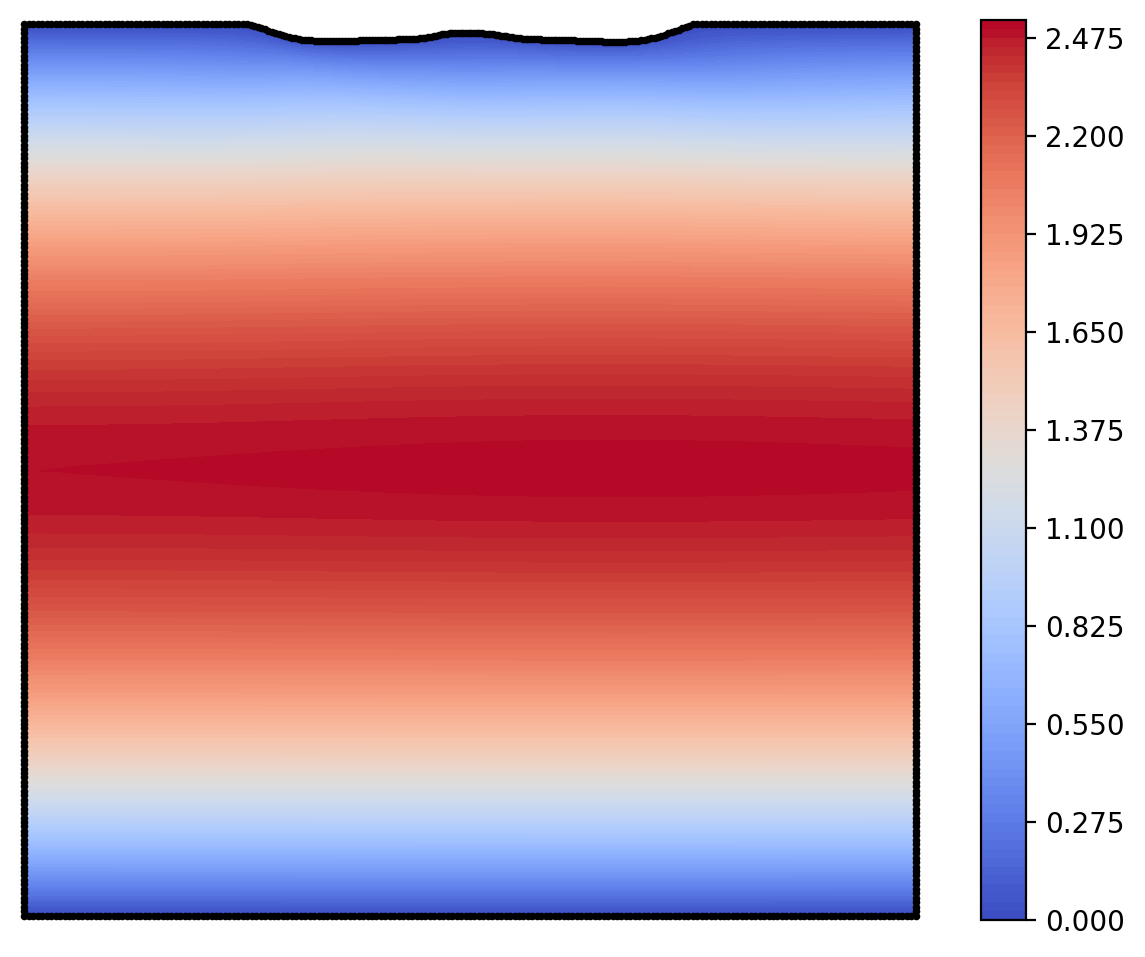}}
  \caption{Channel optimization results by Fireshape and AONN-2. The color represents the velocity magnitude.}
  \label{fig:bump_result_convex}
\end{figure}

As seen from Table~\ref{tab:test4}, AONN-2 is able to achieve lower objective functional values than Fireshape does. Further in Fig.~\ref{fig:bump_result_convex}, we draw the flow fields before and after the optimization for the two initial shapes. Analogous to the previous examples, due to the restriction of  moving mesh, Fireshape usually sacrifices part of the flexibility of the  deformation to maintain the mesh quality (e.g. by adding regularization terms). In comparison, AONN-2 allows the computational domain to be squeezed or stretched more freely, thanks to the considerable flexibility provided by the shape representation points on the free boundary.

\section{Conclusion}
\label{sec:conc}

An adjoint-oriented neural network method for PDE-constrained shape optimization (AONN-2) is proposed in this paper. AONN-2 not only inherits the characteristics from AONN but also takes advantage of shape derivative to optimize the shape represented by the discrete boundary points. Due to the mesh-free nature, AONN-2 can naturally avoid the issues caused by mesh deformation, which are often encountered in mesh-dependent optimization methods. In a series of numerical experiments, we apply AONN-2 to tackle shape optimization problems constrained by the Poisson, Stokes, and Navier-Stokes equations, and compare its performance with the widely-used shape optimization toolbox, Fireshape. The optimization results demonstrate that AONN-2 is able to obtain more desirable results with lower values of objective functionals, and is more robust to various initial shapes.

To further improve the performance of AONN-2, we consider to investigate methods that can reasonably determine the initial shape and position of the optimization object. In addition, we plan to extend the applicability of AONN-2 to topology optimization in the future work.

\section*{CRediT authorship contribution statement}
\textbf{Xili Wang:} Conceptualization, Methodology, Programming, Investigation, Writing – Original draft. \textbf{Pengfei Yin:} Conceptualization, Methodology, Programming, Investigation, Writing - Original draft. \textbf{Bo Zhang:} Conceptualization, Validation, Investigation, Writing – reviewing and editing. \textbf{Chao Yang:} Conceptualization, Validation, Writing – reviewing and editing, Supervision.
\section*{Declaration of competing interest}
The authors declare that they have no known competing financial interests or personal relationships that could have appeared to influence the work reported in this paper.
\section*{Acknowledgments}
This study was supported in part by the National Natural Science Foundation of China (No. 12131002, No. 62306018), and the China Postdoctoral Science Foundation (No. 2022M710211).

\bibliographystyle{model1-num-names}
\bibliography{refs}

\end{document}